\documentclass{siamltex}
\usepackage{graphicx}
\usepackage{mathrsfs,enumerate}
\usepackage{hyperref,amssymb,amsmath,dsfont,eucal,amsfonts}
\newtheorem{remark}[theorem]{Remark}

\newcommand{\nchi}{{\raise.3ex\hbox{$\chi$}}}

\newcommand{\R}{\mathbb{R}}

\newcommand{\N}{\mathbb{N}}

\renewcommand{\aa}{{\mbox{\boldmath$a$}}}

\newcommand{\dD}{{\mbox{\boldmath$D$}}}

\newcommand{\rR}{{\mbox{\boldmath$R$}}}
\newcommand{\mrR}{{\mbox{\boldmath$\mathcal R$}}}

\renewcommand{\tt}{{\mbox{\boldmath$t$}}}

\newcommand{\vv}{{\mbox{\boldmath$v$}}}
\newcommand{\ww}{{\mbox{\boldmath$w$}}}
\newcommand{\wW}{{\mbox{\boldmath$W$}}}
\newcommand{\zz}{{\mbox{\boldmath$z$}}}

\newcommand{\sww}{{\mbox{\scriptsize\boldmath$w$}}}

\newcommand{\tauV}{{\kern-3pt\tau}}

\newcommand{\yy}{{\mbox{\boldmath$y$}}}

\newcommand{\nnu}{{\mbox{\boldmath$\nu$}}}
\newcommand{\snnu}{{\mbox{\scriptsize\boldmath$\nu$}}}

\newcommand{\ssigma}{{\mbox{\boldmath$\sigma$}}}

\newcommand{\zzeta}{{\mbox{\boldmath$ \zeta$}}}
\newcommand{\eeta}{{\mbox{\boldmath$\eta$}}}
\newcommand{\Leb}[1]{{\mathscr L}^{#1}}      

\newcommand{\eps}{\varepsilon}  
\newcommand{\Rd}{{\R^d}}

\newcommand{\Scalar}[2]{\left\langle #1,#2\right\rangle}

\newcommand{\restr}[1]{\lower3pt\hbox{$|_{#1}$}}

\newcommand{\BorelSets}[1]{\mathcal B(#1)}

\newcommand{\Probabilities}[1]{\mathcal P(#1)}          
\newcommand{\ProbabilitiesTwo}[1]{\mathcal P_2(#1)}     

\newcommand{\FinalT}{T}

\newcommand{\weaksto}{{\rightharpoonup^*}}
\newcommand{\weakto}{\rightharpoonup}

\newcommand{\dxi}{{\partial_i}}

\newcommand{\dxj}{{\partial_j}}

\newcommand{\dxij}{{\partial^2_{ij}}}

\newcommand{\tDelta}{\Delta_\gamma}
\newcommand{\tdiv}{\nabla\kern-3pt_\gamma\cdot}

\newcommand{\FPM}{\mathcal M^+}
\newcommand{\RPM}{\mathcal M^+_{\rm loc}}

\newcommand{\FPMMM}[2]{\mathcal P_{#2}}
\newcommand{\CE}{\mathcal{CE}}
\newcommand{\up}{\uparrow}
\newcommand{\down}{\downarrow}
\newcommand{\sfh}{{\sf h}}
\newcommand{\sfm}{{\sf m}}
\newcommand{\sfE}{{\sf E}}

\renewcommand{\d}{\mathrm d}
\newcommand{\D}{\mathrm D}

\newcommand{\email}[1]{E-mail:~\texttt{#1}}
\begin{document}

\title{From Poincar\'e to logarithmic Sobolev inequalities:\\ a gradient flow approach}
\author{Jean Dolbeault\thanks{Ceremade (UMR CNRS no. 7534), Universit\'e Paris-Dauphine, place de Lattre de Tassigny, 75775 Paris C\'edex~16, France. \email{dolbeaul@ceremade.dauphine.fr}}
\and Bruno Nazaret\thanks{Ceremade (UMR CNRS no. 7534), Universit\'e Paris-Dauphine, place de Lattre de Tassigny, 75775 Paris C\'edex~16, France. \email{nazaret@ceremade.dauphine.fr}}
\and Giuseppe Savar\'e\thanks{Universit\`a degli studi di Pavia,
 Dipartimento di Matematica ``F. Casorati'', Via Ferrata 1, 27100, Pavia, Italy. \email{giuseppe.savare@unipv.it}}
}
\date{\today}

\maketitle
\thispagestyle{empty}


\begin{abstract} We use the distances introduced in a previous joint paper to exhibit the gradient flow structure of some drift-diffusion equations for a wide class of entropy functionals. Functional inequalities obtained by the comparison of the entropy with the entropy production functional reflect the contraction properties of the flow. Our approach provides a unified framework for the study of the Kolmogorov-Fokker-Planck (KFP) equation. \end{abstract}

\begin{keywords} 
Optimal transport, Kantorovich-Rubinstein-Wasserstein distance, Generalized Poincar\'e inequality, Continuity equation, Action functional, Gradient flows, Kolmogorov-Fokker-Planck equation
\end{keywords}

\begin{AMS} 
26A51, 26D10, 53D25
\end{AMS}

\section{Setting of the problem}\label{Sec:Intro}

Our starting point concerns nonnegative solutions with finite mass of the heat equation in $\Rd$
\begin{equation}
\label{eq:heat}
\partial_tu_t=\Delta u_t\;.
\end{equation}
It is straightforward to check that for any smooth enough solution of \eqref{eq:heat} and any $C^2$ convex function $\psi$,
\[
\frac d{dt}\int_{\R^d}\psi(u_t)\,dx=-\int_{\R^d}\psi''(u_t)\,|\dD u_t|^2\,dx
\]
so that $\int_{\R^d}\psi(u_t)\,dx$ plays the role of a Lyapunov functional. To extract some information out of such an identity, one needs to analyze the relation between $\int_{\R^d}\psi(u_t)\,dx$ and $\int_{\R^d}\psi''(u_t)\,|\dD u_t|^2\,dx$. This can be done using Green's function or moment estimates, with the drawback that these quantities are explicitly $t$-dependent. It is simpler to rewrite the equation in self-similar variables and replace \eqref{eq:heat} by the Fokker-Planck (FP) equation
\begin{equation}
\label{eq:Fokker-Planck}
\partial_tv_t=\Delta v_t+\nabla\cdot (x\,v)\;.
\end{equation}
This can be done without changing the initial data by the time-dependent change of variables
\[
u_t(x)=\frac 1{R(t)^d}\,v_t\left(\frac x{R(t)}\right),\quad R(t)=\sqrt{1+2t}\;.
\]
We shall restrict our approach to nonnegative initial data $u_0=v_0$. By linearity, we can further assume that
\[
\int_{\R^d}v_t\,dx=\int_{\R^d}u_t\,dx=\int_{\R^d}u_0\,dx=1
\]
without loss of generality. We shall also assume that $\psi$ is defined on $\R^+$. Up to the change of~$\psi$ into $\tilde\psi$ such that $\tilde\psi(s)=\psi(s)-\psi(1)-\psi'(1)(s-1)$, we can also assume that $\psi$ is nonnegative on~$\R^+$ and achieves its minimum value, zero, at $s=1$.

Eq.~\eqref{eq:Fokker-Planck} has a unique nonnegative stationary solution $v=\gamma$ normalized such that \hbox{$\int_{\R^d}\gamma\,dx=1$}, namely
\[
\gamma(x)=\frac{e^{-|x|^2/2}}{(2\pi)^{d/2}}\quad\forall\;x\in\R^d\,.
\]
If we introduce $\rho_t=v_t/\gamma$, then $\rho_t$ is a solution of the Ornstein-Uhlenbeck, or Kolmogorov-Fokker-Planck (KFP), equation
\begin{equation}
\label{eq:Ornstein-Uhlenbeck}
\partial_t\rho_t=\Delta \rho_t-x\cdot\dD \rho_t
\end{equation}
with initial data $\rho_0=v_0/\gamma$. After identifying $\gamma$ with the measure $\gamma\,\Leb d$, the relevant Lyapunov functional, or \emph{entropy}, is $\int_{\R^d}\psi(\rho_t)\,d\gamma$ and
\[
\frac d{dt}\int_{\R^d}\psi(\rho_t)\,d\gamma=-\int_{\R^d}\psi''(\rho_t)\,|\dD \rho_t|^2\,d\gamma\;.
\]
We shall restrict our study to a class of functions $\psi$ for which the entropy and the entropy production functional are related by the inequality
\begin{equation}
\label{Ineq:Entropie-ProductionBeckner}
2\lambda\int_{\R^d}\psi(\rho)\,d\gamma\le\int_{\R^d}\psi''(\rho)\,|\dD \rho|^2\,d\gamma
\end{equation}
for some $\lambda>0$ (it turns out that in the case of the Gaussian measure we can choose $\lambda=1$). This allows us to prove that the entropy is exponentially decaying, namely
\begin{equation}
\label{Eqn:ExponentialDecay}
\int_{\R^d}\psi(\rho_t)\,d\gamma\le\left(\int_{\R^d}\psi(\rho_0)\,d\gamma\right)\;e^{-2\lambda t}\quad\forall\;t\ge0\;,
\end{equation}
if $\rho_t$ is a solution of \eqref{eq:Ornstein-Uhlenbeck} and if $\lambda $ is positive. A sufficient condition for such an inequality is that
\begin{equation}
\label{cdt:concavity}
\mbox{\emph{the function }} h:=1/\psi''\mbox{\emph{ is concave}}
\end{equation}
(see for instance \cite{Arnold-Markowich-Toscani-Unterreiter01}). At first sight, this may look like a technical condition but it has some deep implications. We are indeed interested in exhibiting a gradient flow structure for~\eqref{eq:Fokker-Planck} associated with the entropy or, to be more precise, to establish that, for some distance, the gradient flow of the entropy is actually \eqref{eq:Fokker-Planck}. It turns out that \eqref{cdt:concavity} is the natural condition as we shall see in Section \ref{subsec:ActionDensity}.

\medskip The entropy decays exponentially according to \eqref{Eqn:ExponentialDecay} not only when one considers the $L^2_\gamma(\R^d)$ norm (the norm of the square integrable functions with respect to the Gaussian measure $\gamma$), i.e.~the case $\psi(\rho)=(\rho-1)^2/2$, or the classical entropy built on $\psi(\rho)=\rho\,\log\rho$, for which \eqref{eq:Ornstein-Uhlenbeck} is the gradient flow with respect to the usual Wasserstein distance (according to the seminal paper \cite{Jordan-Kinderlehrer-Otto98} of Jordan, Kinderlehrer and Otto). We also have an exponential decay result of any entropy generated by
\[
\psi(\rho)=\frac{\rho^{2-\alpha}-1-(2-\alpha)(\rho-1)}{(2-\alpha)(1-\alpha)}=:\psi_\alpha(\rho)\;,\quad\alpha\in[0,1)\;,
\]
and more generally any $\psi$ satisfying \eqref{cdt:concavity}. Notice by the way that $\psi(\rho)=\psi_\alpha(\rho)$ is compatible with \eqref{cdt:concavity} if and only if $\alpha\in[0,1)$ and that $\psi(\rho)=\rho\,\log\rho$ appears as the limit case when $\alpha\to 1_-$.

The exponential decay is a striking property which raises the issue of the hidden mathematical structure, a question asked long ago by F. Poupaud. As already mentionned, the answer lies in the gradient flow interpretation and the construction of the appropriate distances. Such distances, based on an action functional related to~$\psi$, have been studied in \cite{MR2448650}. Our purpose is to exploit this action functional for the construction of gradient flows, not only in the case corresponding to \eqref{eq:Ornstein-Uhlenbeck} but also for KFP equations based on general $\lambda$-convex potentials $V$. For the convenience of the reader, the main steps of the strategy have been collected in Section~\ref {Sec:Formal}, without technical details (for instance on the measure theoretic aspects of our approach).

Coming back to our basic example, namely the solution of \eqref{eq:Ornstein-Uhlenbeck}, we may observe that a solution can easily be represented using the Green kernel of the heat equation and our time-dependent change of variables. If $\psi(\rho)=\psi_\alpha(\rho)$, $\alpha\in[0,1)$, we may observe that the exponential decay of the entropy can be obtained using the known properties of the heat flow and the homogeneity of $\psi_\alpha$, while the contraction properties of the heat flow measured in the framework of the weighted Wasserstein distances introduced in \cite{MR2448650} can be translated into the exponential decay of the distance of the solution of \eqref{eq:Ornstein-Uhlenbeck} to the gaussian measure $\gamma$, if we assume that $\rho\,\gamma$ is a probability measure. We shall however not pursue in this direction as it is very specific of the potential $V(x)=\frac 12\,|x|^2$ and of the heat flow (for which an explicit Green function is available).

\medskip Let us conclude this introductory section by a brief review of the literature on the functional inequalities based on entropies such that \eqref{cdt:concavity} holds. Such functionals are sometimes called \emph{$\varphi$-entropies}. In this paper, we shall however avoid this denomination to prevent from possible confusions with the function $\phi$ and the functional $\Phi$ used below to define the action and the weighted Wasserstein distances~$W_h$. 

We shall refer to \cite{MR2081075,MR1796718} for a probabilistic point of view. A proof of \eqref{Eqn:ExponentialDecay} under Assumption \eqref{cdt:concavity} and an hypothesis of convexity of $V$ can be found for instance in \cite{Arnold-Markowich-Toscani-Unterreiter01} or in the more recent paper \cite{MR2609029}. This approach is based on the Bakry-Emery method \cite{MR889476,MR2273884} and heavily relies on the flow of KFP or, equivalently, on the geometric properties of the Ornstein-Uhlenbeck operator (using the \emph{carr\'e du champ:} see \cite{MR2609029}). Strict convexity of the potential is usually required, but can be removed afterwards by various methods: see \cite{Arnold-Markowich-Toscani-Unterreiter01,MR2317340,MR2435196}. For capacity-measure approaches of \eqref{Ineq:Entropie-ProductionBeckner}, we shall refer to \cite{MR2320410,MR2346509,MR2366398}. The inequality \eqref{Ineq:Entropie-ProductionBeckner} itself has been introduced in \cite{Beckner89} with a proof based on the hypercontractivity of the heat flow and spectral estimates, and later refined and adapted to general potentials in \cite{MR2375056}.

Concerning gradient flows and distances of Wasserstein type, there has been a huge activity over the last years. We can refer to \cite{Jordan-Kinderlehrer-Otto98,Benamou-Brenier00} for fundamental ideas, and to two books, \cite{ags,MR2459454}, for a large overview of the field. Many other contributions in this area will be quoted whenever needed in the proofs.

\section{Formal point of view: definitions, strategy and main results}\label{Sec:Formal}

In Section~\ref{Sec:Intro}, we have considered the case of the harmonic potential $V(x)=\frac 12\,|x|^2$. We generalize the setting to any smooth, convex \emph{potential} $V:\Rd\to\R$ with
\begin{equation}
\label{eq:Hessian}
D^2 V\ge \lambda\,\mathsf I\;,\quad \lambda\ge0\;,
\end{equation}
and consider the \emph{reference measure} $\gamma$ given by
\begin{equation}
\label{eq:DensityLebesgue}
\gamma:=e^{-V}\,\Leb d
\end{equation}
where $\Leb d$ denotes Lebesgues's measure on $\Rd$. We assume that
\begin{equation}
\label{eq:PartitionFunctionFinite}
\gamma(\Rd)=\int_{\Rd}e^{-V}\,dx=:Z<\infty\;.
\end{equation}
Next we define the \emph{action density} $\phi:(0,\infty)\times \Rd\to\R$ as
\[
\label{eq:47}
\phi(\rho,\ww):=g(\rho)\,|\ww|^2=\frac{|\ww|^2}{h(\rho)}
\]
for some concave, positive, non decreasing function $h$ with \emph{sublinear} growth. The function $g$ is therefore convex and also satisfies the condition
\begin{equation}
\label{eq:2}
2\,(g')^2\le g\,g''\,.
\end{equation}
Our main example is $h(\rho):=\rho^\alpha$ for some $\alpha\in(0,1)$.
Based on the action density, we can define the \emph{action functional} by
\begin{equation}
\label{eq:53}
\Phi(\rho,\ww):=\int_{\Rd}\phi(\rho,\ww)\,d\gamma\;.
\end{equation}

\medskip\noindent{\bf The Kolmogorov-Fokker-Planck (KFP) equation.}
With the notations $\tDelta:=\Delta-\dD V\cdot \dD$, the equation
\begin{equation}
\label{eq:51}
\partial_t\rho_t-\tDelta\rho_t=0
\end{equation}
determines the \emph{Kolmogorov-Fokker-Planck (KFP) flow} $S_t:\rho_0\mapsto \rho_t$. Its first variation, $\rR_t:\ww_0\mapsto \ww_t$, can be obtained as the solution of the \emph{modified Kolmogorov-Fokker-Planck equation}
\[
\label{eq:52}
\partial_t\ww_t-\tDelta \ww_t+D^2 V\,\ww_t=0\;.
\]
If $\ww_0=\dD\rho_0$, then $\ww_t=\dD\rho_t$, which can be summarized by
\[
\label{eq:55}
\dD (S_t\rho_0)=\rR_t(\dD\rho_0)\;.
\]
By duality, using the notations $\tdiv \ww:=\nabla\cdot\ww-\dD V\cdot\ww$ and $\nabla\cdot\ww:=\sum_{i=1}^d\partial\ww_i/\partial x_i$, if $\tdiv \ww_0=\rho_0$, we also find that $\tdiv\ww_t=\rho_t$, which amounts to
\begin{equation}
\label{eq:56}
\tdiv(\rR_t\ww_0)=S_t(\tdiv \ww_0)
\end{equation}
(see Theorem~\ref{thm:density_gradient_link} for details). If $\mu=\rho\,\gamma$, we define the semigroup $\mathcal S_t$ acting on measures by $\mathcal S_t\mu:=(S_t\rho)\,\gamma$.

Consider an \emph{entropy density} function $\psi$ such that $\psi(1)=\psi'(1)=0$. If we define the \emph{entropy} functional by
\[
\label{eq:60}
\Psi(\rho):=\int_{\Rd}\psi(\rho)\,d\gamma
\]
and the \emph{entropy production}, or generalized Fisher information functional, as the action functional for the particular choice $\ww=\dD\rho$, i.e.
\[
\label{eq:57}
P(\rho):=\Phi(\rho,\dD\rho)\;,
\]
then, along the KFP flow, we get
\begin{equation}
\label{eq:58}
\frac d{dt}\Psi(\rho_t)=-P(\rho_t)=-\Phi(\rho_t,\dD\rho_t)
\end{equation}
for a solution $\rho_t$ of \eqref{eq:51} if
\[
\label{eq:59}
\psi''=g\;.
\]
Notice that \eqref{cdt:concavity} and \eqref{eq:2} are equivalent. See Section~\ref {subsec:ActionDensity} for more details. The main estimate for this paper goes as follows.
\smallskip\begin{theorem}
\label{Thm:MainEstimate}
Under Assumptions \eqref{eq:Hessian}--\eqref{eq:2}, if $\Phi(\rho_0,\ww_0)<\infty$, $\rho_t=S_t\rho_0$ and $\ww_t=\rR_t\ww_0$, then
\[
\label{eq:54}
\frac d{dt}\,\Phi(\rho_t,\ww_t)+2\lambda\,\Phi(\rho_t,\ww_t)\le0\quad\forall\;t\ge 0\;.
\]
\end{theorem}\smallskip
In particular the action functional decays exponentially if $\lambda$ is positive:
\begin{equation}
\label{eq:62}
\Phi(\rho_t,\ww_t)\le e^{-2\lambda t}\,\Phi(\rho_0,\ww_0)\quad\forall\;t\ge0\;.
\end{equation}
At formal level, this follows by an easy convexity argument. The rigorous proof requires many regularizations. See Theorem~\ref{thm:action_decay} for a more detailed version of this result. Now let us review some of the consequences of Theorem~\ref{Thm:MainEstimate}.

\medskip\noindent{\bf Entropy, entropy production and generalized Poincar\'e inequalities.} We can now apply Theorem~\ref{Thm:MainEstimate} to the KFP flow. With $\ww=\dD\rho$, we find that the \emph{entropy production} functional decays exponentially:
\begin{equation}
\label{eq:63}
\frac d{dt}P(\rho_t)+2\lambda\,P(\rho_t)\le 0\;,\quad P(\rho_t)\le e^{2\lambda t} P(\rho_0)\quad\forall\;t\ge 0
\end{equation}
if $\lambda$ is positive. By integrating \eqref{eq:58} along the KFP flow when $t$ varies in $\R^+$, using \eqref{eq:63} and $\Psi(1)=0$, we recover for $\rho=\rho_0$ the \emph{generalized Poincar\'e inequalities}
\begin{equation}
\label{eq:64}
\Psi(\rho)\le \frac 1{2\lambda}\,P(\rho)
\end{equation}
found by Beckner in \cite{Beckner89} in the case of the harmonic potential and for $h(\rho):=\rho^\alpha$, $\alpha\in(0,1)$, and generalized for instance in \cite{Arnold-Markowich-Toscani-Unterreiter01}. Such inequalities interpolate between Poincar\'e and logarithmic Sobolev inequalities.

If we combine \eqref{eq:64} with \eqref{eq:63}, we find that the \emph{entropy} decays according to
\[
\label{eq:65}
\frac d{dt}\Psi(\rho_t)+2\lambda\Psi(\rho_t)\le 0\;,\quad\Psi(\rho_t)\le e^{-2\lambda t}\,\Psi(\rho_0)\quad\forall\;t\ge 0\;.
\]
By integrating from $0$ to $t$ the inequality
\[
\label{eq:67}
\frac{d}{dt}\Big(t\,P(\rho_t)\Big)=P(\rho_t)+t\,\frac d{dt}P(\rho_t)\le P(\rho_t)=-\frac d{dt}\Psi(\rho_t)\;,
\]
which itself follows from \eqref{eq:58} and \eqref{eq:63}, we observe a \emph{first regularization effect} along the KFP flow, namely
\begin{equation}
\label{eq:66}
t\,P(\rho_t)\le \Psi(\rho_0)\quad\forall\;t\ge 0\;.
\end{equation}
If $\lambda$ is positive, we can refine this estimate and actually prove by the same method that $\frac{e^{2\lambda t}-1}{2\lambda}\,P(\rho_t)\le \Psi(\rho_0)$ for any $t\ge 0$.

\medskip\noindent{\bf The $h$-Wasserstein distance.} If $\mu$ is a measure with absolutely continuous part~$\rho$ with respect to $\gamma$, and singular part $\mu^\perp$, if $\nnu$ is a vector valued measure which is absolutely continuous with respect to $\gamma$ and has a modulus of continuity $\ww$, i.e.~if
\begin{equation}
\label{eq:decomp}
\mu=\rho\,\gamma+\mu^\perp\quad\mbox{and}\quad\nnu=\ww\,\gamma\;,
\end{equation}
we can extend the action functional $\Phi$ to the measures $\mu$ and $\nnu$ by setting
\[
\label{eq:48}
\Phi(\mu,\nnu)=\Phi(\rho,\ww)=\int_{\Rd}\phi(\rho,\ww)\,d\gamma\;.
\]
We shall say that there is an \emph{admissible path} connecting $\mu_0$ to $\mu_1$ if there is a solution $(\mu_s,\nnu_s)_{s\in [0,1]}$ to the \emph{continuity equation}
\[
\label{eq:50}
\partial_s \mu_s+\nabla\cdot\nnu_s=0\;,\quad s\in [0,1]\;,
\]
and will denote by $\Gamma(\mu_0,\mu_1)$ the set of all admissible paths. With these tools, we can define the \emph{$h$-Wasserstein distance} between $\mu_0$ and $\mu_1$ by
\[
\label{eq:49}
W^2_h(\mu_0,\mu_1):=\inf\Big\{\int_0^1 \Phi(\mu_s,\nnu_s)\,ds\,:\,(\mu,\nnu)\in\Gamma(\mu_0,\mu_1)\Big\}\,.
\]
Notice that $h$ in ``$h$-Wasserstein distance'' refers to the dependence of $\Phi$ in $h$ through the action density $\phi$, the usual Wasserstein distance corresponding to $h(\rho)=\rho$. If $(\mu_t)_{t\in (0,T)}$ is a curve of measures, its \emph{$h$-Wasserstein velocity} $|\dot\mu_t|$ is determined~by
\[
|\dot\mu_t|^2=\inf_\nnu\Big\{\Phi(\mu,\nnu)\,:\,\nabla\cdot \nnu=-\,\partial_t\mu_t\Big\}\,.
\]
Using the decomposition \eqref{eq:decomp}, we compute the \emph{derivative of the entropy along the curve} $(\mu_t)_{t\in (0,T)}$ as
\begin{displaymath}
\frac d{dt}\Psi(\rho_t)=\int_{\Rd} \psi'(\rho_t)\,\partial_t\rho_t\,d\gamma=\int_{\Rd} \psi''(\rho_t)\,\dD\rho_t\cdot \ww_t\,d\gamma
\end{displaymath}
and find that
\begin{equation}
\label{eq:69bis}
-\frac d{dt}\Psi(\rho_t)=-\int_{\Rd} \sqrt{\psi''(\rho_t)}\,\dD\rho_t\cdot \sqrt{\psi''(\rho_t)}\,\ww_t\,d\gamma\le\sqrt {P(\rho_t)}\,|\dot\mu_t|
\end{equation}
by the Cauchy-Schwarz inequality. Along the KFP flow, we know that
\[
\label{eq:70}
\frac d{dt}\Psi(\rho_t)=-P(\rho_t)=-|\dot \mu_t|^2=-\sqrt {P(\rho_t)}\,|\dot\mu_t|\;,
\]
which is the equality case in \eqref{eq:69bis}. This characterizes the KFP flow as the steepest descent flow of the entropy $\Psi$, i.e.~this is a first charaterization of KFP as \emph{the gradient flow of $\Psi$} with respect to the $h$-Wasserstein distance.

The KFP flow connects $\mu=\rho\,\gamma$ with $\mu_\infty=\gamma$ and it has been established in \cite{MR2448650} that one can estimate the length of the path by
\begin{equation}
\label{eq:72}
W_h(\mu,\gamma)=\int_0^{\infty}\sqrt{P(\rho_t)}\,dt=\int_0^{\infty}|\dot \mu_t|\,dt
\end{equation}
(see Section~\ref{subsec:WeightedWasserstein} for details). According to \eqref{eq:63}, we get
\[
\label{eq:72bis}
W_h(\mu,\gamma)\le\sqrt {P(\rho)}\int_0^{\infty}e^{-\lambda t}\,dt=\frac 1\lambda\,\sqrt{P(\rho)}\;.
\]
This establishes the \emph{entropy production -- distance estimate}
\[
\label{eq:71}
W_h(\mu,\gamma)\le \frac 1\lambda\,\sqrt{P(\rho)}\;,\quad\mbox{if}\quad\mu=\rho\,\gamma\;.
\]
Along the KFP flow, we also find that
\[
\label{eq:74}
-\frac{d}{dt}\sqrt{\Psi(\rho_t)}=\frac{P(\rho_t)} {2\sqrt{\Psi(\rho_t)}}\ge\sqrt {\frac\lambda2\,P(\rho_t)}
\]
using \eqref{eq:64}. By applying \eqref{eq:72}, this establishes the \emph{(Talagrand) entropy -- distance estimate}
\[
\label{eq:73}
W^2_h(\mu,\gamma)\le \frac 2\lambda\,\Psi(\rho)\;.
\]

\medskip\noindent{\bf Contraction properties and gradient flow structure.} Here as in \cite{MR2448650}, we use the technique introduced in \cite{MR2192294} and extended in \cite[\S\,2]{MR2452882}: we consider a geodesic (or an approximation of a geodesic), and evaluate the derivative of the action functional along a family of curves obtained by evolving the geodesic with the KFP flow.

Consider an $\varepsilon$-geodesic $(\rho^s,\ww^s)$ connecting $\mu^0=\rho^0\,\gamma$ to $\mu^1=\rho^1\,\gamma$, i.e.~an admissible path in $\Gamma(\mu_0,\mu_1)$ such that $\Phi(\rho^s_0,\ww^s_0)\le W_h^2(\rho^0_0,\rho_0^1)+\varepsilon$ for any $s\in(0,1)$ and observe that by~\eqref{eq:56}, we know that $(\rho^s_t=S_t\rho^s,\ww^s_t=\rR_t\ww^s)$ is still an admissible curve connecting $S_t\rho^0$ to $S_t\rho^1$. Therefore \eqref{eq:62} yields
\[
\label{eq:77}
W_h^2(\rho^0_t,\rho_t^1)\le\int_0^1\Phi(\rho^s_t,\ww^s_t)\,ds\le e^{-2\lambda t}\int_0^1
\Phi(\rho^s_0,\ww^s_0)\,ds\le e^{-2\lambda t}\left(W_h^2(\rho^0_0,\rho_0^1)+\varepsilon\right)\;,
\]
which, by letting $\varepsilon\to0$, proves that the KFP flow \emph{contracts the distance:}
\[
\label{eq:76}
W_h(\mathcal S_t\mu^0,\mathcal S_t\mu^1)\le e^{-\lambda t}\,W_h(\mu^0,\mu^1)\quad\forall\;t\ge 0\;.
\]
See Theorem~\ref{thm:main_contraction} for more details.

Next, we should again consider an $\varepsilon$-geodesic, but for simplicity we assume that there is a geodesic $(\rho^s,\ww^s)$ connecting $\sigma=\mu^0=\rho^0\,\gamma$ to $\mu=\mu^1=\rho^1\,\gamma$, i.e.~such that $\Phi(\rho^s,\ww^s)=W_h^2(\sigma,\mu)$, and consider the path
\[
\label{eq:78}
(\rho_t^s,\ww_t^s):=(S_{st}\rho^s,\rR_{st}\ww_s+t\,\dD\rho_t^s)
\]
connecting $\sigma$ to $\mu_t:=\mathcal S_t\mu$. Notice that our notations mean that $\rho^s=\rho^s_0$. Since
\[
\label{eq:79}
\partial_s\rho_t^s=\rho_t^s+t\,\tDelta \rho_t^s=\tdiv (\ww_t^s+t\,\dD\rho_t^s)\;,
\]
the path is admissible and, as a consequence,
\[
\label{eq:80}
W_h^2(\mu_t,\sigma)\le\int_0^1\Phi(\rho_t^s,\ww_t^s)\,ds\;.
\]
We can therefore differentiate the right hand side in the above inequality instead of the distance and furthermore notice that it is sufficient to do it at $t=0$; see Theorem~\ref{thm:metric_characterization} and its proof for details. Along the KFP flow we find that
\begin{equation}
\label{eq:75}
\frac 12\,\frac d{dt}W_h^2(\mu_t,\sigma)+\frac \lambda2\,W_h^2(\mu_t,\sigma) \le \Psi(\sigma\,|\,\gamma)-\Psi(\mu_t\,|\,\gamma)\;.
\end{equation}
This is the strongest metric formulation of a \emph{$\lambda$-contracting gradient flow}. Here we have defined the relative entropy as $\Psi(\mu\,|\,\gamma):=\psi(\rho)$ if $\mu\ll\gamma$ and $\mu=\rho\,\gamma$, and $\Psi(\sigma\,|\,\gamma):=+\infty$ otherwise. Hence we recover a second characterization of the fact that KFP \emph{is the gradient flow of $\Psi$} with respect to $W_h$.

As another consequence, the entropy \emph{$\Psi$ is geodesically $\lambda$-convex}. This follows from \eqref{eq:75}. Fix a geodesic $\mu^s$ between $\mu^0$ and $\mu^1$, follow the evolution of $\mu^s$ by KFP taking first $\mu^0$ and then $\mu^1$ fixed, and apply \eqref{eq:75} with $\mu_t:=\mathcal S_t\mu^s$ and $\mu=\mu^0$ or $\mu=\mu^1$. Because of the minimality of the energy along the geodesic at time $t=0$, by summing the two resulting inequalities we prove the convexity inequality of $\Psi$. See \cite[Theorem 3.2]{MR2452882} for more details.

As a final observation, let us notice that, directly from the metric formulation \eqref{eq:75}, it follows that the KFP flow also has the following \emph{regularizing properties:}
\[
\label{eq:81-82}
\Psi(\rho_t)\le \frac 1{2t}\,W_h^2(\rho_0,\gamma)\quad\mbox{and}\quad P(\rho_t)\le \frac 1{t^2}\,W_h^2(\rho_0,\gamma)\quad\forall\;t\ge 0\;.
\]
The first estimate can indeed be obtained by
integrating \eqref{eq:75} (with $\lambda=0$ and $\sigma=\gamma$)
from $0$ to $t$ and recalling that $t\mapsto \Psi(\rho_t)$ is
decreasing.
As for the second one, we observe that also $t\mapsto P(\rho_t)$ is
decreasing by \eqref{eq:63}, so that \eqref{eq:66} and \eqref{eq:75} yield
\begin{displaymath}
 \frac d{dt}\left(\frac {t^2}2\,P(\rho_t)\right)\le t\,P(\rho_t) \le
 \Psi(\rho_t)\le -\frac 12\,\frac d{dt}W_h^2(\mu_t,\gamma)\;.
\end{displaymath}
A further integration in time from $0$ to $t$ completes the proof.
Notice that it is crucial to start from a measure $\mu=\rho_0\,\gamma$ at \emph{finite} distance from~$\gamma$.

\section {Definition and properties of the weighted Wasserstein distance}\label{Sec:Known}

In this section we first recall some definitions and results taken from \cite{MR2448650}. The measure~$\gamma$ and the functions $\phi$ and $\psi$ are as in Section~\ref{Sec:Formal}, and we assume that Conditions~\eqref{eq:Hessian}--\eqref{eq:2} are satisfied.

\subsection{Properties of the potential}

Let $V:\Rd\to\R$ be a $\lambda$-convex and continuous
potential. We assume that $\lambda$ is nonnegative and $\lambda$-convexity means that the map $x\mapsto V(x)-\frac\lambda2\,|x|^2$ is convex. When $V$ is smooth in $\Rd$, this condition is equivalent to \eqref{eq:Hessian}. We are assuming that $e^{-V}$ is integrable in $\Rd$, so that we can introduce the finite, positive, log-concave measure $\gamma$ defined by \eqref{eq:DensityLebesgue}. For simplicity, we shall assume that $\gamma$ is a probability measure, i.e.~$Z=1$, which can always be enforced by replacing $V$ by $V+\log Z$. The potential $V$ being convex, the integrability of $e^{-V}$ is equivalent to the property that $V(x)\uparrow\infty$ at least linearly as $|x|\uparrow\infty$; see e.g.~\cite[Appendix]{Ambrosio-Savare-Zambotti07}. As a consequence, there exist two constants $A>0$, $B\ge 0$ such that
\begin{equation}
\label{eq:scheme3:3}
V(x)\ge A\,|x|-B\quad\forall\;x\in\Rd\,.
\end{equation}

We recall that non smooth, convex potentials $V$ can be approximated from below by an increasing sequence of convex potentials $V_n$:
\[
\label{eq:scheme3:27}
V_n(x):=\frac \lambda2\,|x|^2+\inf_{y\in\Rd} \Big(\frac n2\,|x-y|^2+V(y)-\frac \lambda2\,|y|^2\Big)\,.
\]
Moreover, the potentials $V_n$ are $\lambda$-convex and, even in the
case $\lambda=0$, they satisfy conditions \eqref{eq:scheme3:3} with
respect to constants $A$ and $B$ which are \emph{independent of $n$.}
In particular,
the log-concave measures $\gamma_n:=e^{-V_n}\Leb d$ weakly$^*$ and
monotonically converge in $C^0_b(\Rd)'$ to~$\gamma$.
By this regularization techniques, many results could be
extended to the case when $V$ is just lower semicontinuous and can
take the value $+\infty$.

\subsection{Convexity of the action density}\label{subsec:ActionDensity}
As in Section~\ref{Sec:Formal}, consider $g$ and $h$ on $(0,\infty)$ such that $g(\rho)=1/h(\rho)$ and $\phi(\rho,\ww)=g(\rho)\,|\ww|^2=|\ww|^2/h(\rho)$. The following result has already been observed in \cite{MR2448650} but we reproduce it here for completeness.
\smallskip\begin{lemma}[Convexity of the action density]
\label{le:action1}
With the notations of Section~\ref{Sec:Formal}, the action density $\phi$ is \emph{convex} if and only if $h$ is concave on $(0,\infty)$ or, equivalently, if $g$ satisfies Condition~\eqref{eq:2}.
\end{lemma}\smallskip
\begin{proof} By standard approximations, it is not restrictive to assume that $g$, $h\in C^2(0,\infty)$. First of all observe that
\[
g^3\,h''=2\,(g')^2-g\,g''\,,
\]
so that $h''$ is nonpositive if and only if $2\,(g')^2\le g\,g''$. Next we evaluate the second derivative of $\phi$ along the direction of the vector $\zz=(x,\yy)\in\R\times\Rd$ as
\[
\Scalar{\D^2\phi(\rho,\ww)\,\zz}\zz=g''(\rho)\,|\ww|^2\,x^2+4\,g'(\rho)\,\ww\cdot x\,\yy+2\,g(\rho)\,|\yy|^2\,.
\]
By minimizing with respect to $x\in\R$, we get
\begin{equation}
\label{eq:Cap2:32}
g''(\rho)\,|\ww|^2\,\Scalar{\D^2\phi(\rho,\ww)\,\zz}\zz\ge2\,\left[g''(\rho)\,|\ww|^2\,g(\rho)\,|\yy|^2-2\,(g'(\rho)\,\ww\cdot\yy)^2\right]
\end{equation}
if $g''(\rho)>0$, with equality for the appropriate choice of $x$. The convexity of $\phi$ is thus equivalent to
\[
\label{eq:Cap2:29}
g''(\rho)\,|\ww|^2\,g(\rho)\,|\yy|^2\ge2\,(g'(\rho)\,\ww\cdot\yy)^2\quad\forall\;\rho>0\,,\quad\forall\;\yy\,,\;\ww\in\Rd\,.
\]
If $\phi$ is convex, by choosing $\yy:=h(\rho)\,g'(\rho)\,\ww$ and using $h(\rho)\,g(\rho)=1$, we get
\[
\label{eq:Cap2:29bis}
g''(\rho)\,|\ww|^2\,h(\rho)\,(g'(\rho))^2\,|\ww|^2\ge2\,\left[h(\rho)\,(g'(\rho))^2\,|\ww|^2\right]^2\quad\forall\;\rho>0\,,\quad\forall\;\ww\in\Rd\,,
\]
which yields \eqref{eq:2}. Conversely, the convexity of $\phi$ follows from $(\ww\cdot\yy)^2\le|\ww|^2\,|\yy|^2$.
\end{proof}

\medskip We can introduce a \emph{modulus of convexity} as follows. Assume that for some $\alpha\in(0,1]$ we have
\begin{equation}
\label{eq:Cap2:31}
g(\rho)\,g''(\rho)\ge(1+\alpha^{-1})\,(g'(\rho))^2\quad\forall\;\rho>0\;.
\end{equation}
By \eqref{eq:Cap2:32}, we obtain the refined estimate
\begin{equation}
\label{eq:scheme:14}
\Scalar{\D^2\phi(\rho,\ww)\,\zz}\zz\ge 2\,\beta\,\phi(\rho,\yy)\quad\forall\;\zz=(x,\yy)\in\R^{d+1}\,,\quad\mbox{with}\;\beta:=\frac{1-\alpha}{1+\alpha}\;.
\end{equation}
Such a refinement has interesting consequences, which have been investigated in \cite{MR2375056,MR2152502,MR2435196}. The refined convexity assumption \eqref{eq:Cap2:31} is equivalent to
\[
\label{eq:cap3:3}
h^{1/\alpha}\;\mbox{is concave}\,.
\]
\begin{remark}[Main example]
\label{rem:main_example}
Our main example is provided by the function
\[
\label{eq:scheme3:1}
h(\rho):=\rho^\alpha\,,\quad 0\le \alpha\le 1\;,\quad\phi(\rho,\ww)=\frac{|\ww|^2}{\rho^\alpha}\;,
\]
which satisfies \eqref{eq:scheme:14}. When $\alpha=0$ we simply get
\[
\label{eq:Cap2:24}
\phi(\rho,\ww):=|\ww|^2\,,
\]
and for $\alpha=1$ we have the $1$-homogeneous functional
\[
\label{eq:Cap2:25}
\phi(\rho,\ww):=\frac{|\ww|^2}\rho\;.
\]
\end{remark}
Notice that the above considerations can be generalized to matrix-valued functions $g$ and $h$: see \cite[Example 3.4]{MR2448650}.

\subsection{The action functional on densities}

The \emph{action functional} $\Phi$ induced by $\phi$ has been defined by \eqref{eq:53}, with domain
\[
\label{eq:Cap2:9}
\mathcal D(\Phi):=\Big\{(\rho,\ww)\in L^1_\gamma(\Rd)\times L^1_\gamma(\Rd;\Rd)\,:\,\rho\ge0\,,\;\Phi(\rho,\ww)<\infty\Big\}\,.
\]
Assuming as in Section~\ref{subsec:ActionDensity} that $\phi$ convex, it is well known that if $(\rho_k)_{k\in\N}$ and $(\ww_k)_{k\in\N}$ are such that $(\rho_k,\ww_k)\in\mathcal D(\Phi)$ for any $k\in\N$ and if $\rho_k\weakto\rho$ in $L^1_\gamma(\Rd)$, and $\ww_k\weakto^* \ww\in L^1_\gamma(\Rd;\Rd)$ as $n\uparrow\infty$, then by lower semi-continuity of $\Phi$, we have
\[
\label{eq:Cap2:7}
\liminf_{n\uparrow \infty}\Phi(\rho_k,\ww_k)\ge\Phi(\rho,\ww)\;.
\]
\smallskip\begin{lemma}[Approximation by smooth bounded densities]
\label{le:smooth_approximations}
Consider two functions $\rho\in L^1_\gamma(\Rd)$ and $\ww\in L^1_\gamma(\Rd;\Rd)$ such that $\rho\ge0$ and $\Phi(\rho,\ww)<\infty$. Then there exist two sequences $(\rho_k)_{k\in\N}$ and $(\ww_k)_{k\in\N}$ of bounded smooth functions (with bounded derivatives of arbitrary orders) such that $\inf_\Rd \rho_k>0$ and
\[
\label{eq:scheme3:13}
\begin{gathered}
\lim_{k\up\infty}\rho_k=\rho\quad\mbox{in}\;L^1_\gamma(\Rd)\;,\quad\lim_{k\up\infty}\ww_k=\ww\quad\mbox{in}\;L^1_\gamma(\Rd;\Rd)\;,\\
\int_\Rd\rho_k\,d\gamma=\int_\Rd\rho\,d\gamma\quad\forall\;k\in\N\quad\mbox{and}\quad\lim_{k\up\infty}\int_\Rd\phi(\rho_k,\ww_k)\,d\gamma= \int_\Rd \phi(\rho,\ww)\,d\gamma\;.
\end{gathered}
\]
\end{lemma}\smallskip
\begin{proof} We first truncate $\rho$ and $\ww$ from above as follows. Let $m:=\int_\Rd \rho\,d\gamma$ and, for any $k\in\N$, $ m_k:=\int_\Rd (\rho\land k)\,d\gamma,$ $R_k:=\{x\in\Rd\,:\,\rho(x)\le k\}$. We set $\rho_k:=m_k^{-1}\,m\,(\rho\land k)$ and
\[
w_k(x):=\begin{cases}
w(x)&\mbox{if}\;|w(x)|\le k\mbox{ and}\;x\in R_k\;,\\
0&\mbox{otherwise}\;.
\end{cases}
\]
Clearly $\rho_k\to\rho$, $\ww_k\to\ww$ pointwise $\gamma$ a.e.~in $\Rd$, so that Fatou's Lemma yields
\begin{equation}
\label{eq:scheme3:15}
\liminf_{k\uparrow\infty}\Phi (\rho_k,\ww_k)=\Phi(\rho,\ww)\;.
\end{equation}
Since $\rho\land k\to\rho$ in $L^1_\gamma(\Rd)$ as $k\up\infty$, we have $m_k\to m$ and $\rho_k\to\rho$ in $L^1_\gamma(\Rd)$. The dominated convergence theorem also yields $\ww_k\to\ww$ in $L^1_\gamma(\Rd;\Rd)$. Finally, since $\rho_k\ge\rho$ and $|\ww_k|\le |\ww|$ on~$R_k$, and since $g$ is non increasing,
\begin{multline*}
\label{eq:scheme3:14}
\Phi(\rho_k,\ww_k)=\!\int_\Rd \phi(\rho_k,\ww_k)\,d\gamma=\!\int_{R_k}\phi(\rho_k,\ww_k)\,d\gamma\\
\le\!\int_{R_k}\phi(\rho,\ww)\,d\gamma\le \!\int_{\Rd}\phi(\rho,\ww)\,d\gamma=\Phi(\rho,\ww)\;,
\end{multline*}
so that the ``$\liminf$'' in \eqref{eq:scheme3:15} is in fact a limit.

Next we perform a lower truncation on $\rho$. By a diagonal argument, it is sufficient to approximate the functions $\rho_k$ and $\ww_k$ we have just introduced, so we can assume that $\rho$ is essentially bounded by a constant $k$ and we omit the dependence on $k$. For $\delta>0$ we now set $\rho_\delta:=(\rho+\delta)\,m/(m+\delta)$. Observe that
\[
\rho_\delta-m=\frac m{m+\delta}\,(\rho-m)\quad\mbox{and}\quad\rho-\rho_\delta=\frac{\delta}{m+\delta}\,(\rho-m)
\]
so that $m\le\rho_\delta\le\rho$ on the set $R_m^c$ and, by convexity of $g$, we get
\[
g(\rho_\delta)\le C_\delta\,g(\rho)\quad\mbox{where}\quad C_\delta=1+\delta\,\frac{|g'(m)|\,(k-m)}{g(k)\,(\delta+m)}\;.
\]
On the other hand, on the set $R_m$, we have $\rho\le\rho_\delta$, and then $g(\rho_\delta)\le g(\rho)$. As a consequence,
\[
\int_\Rd\phi(\rho_\delta,\ww)\,d\gamma\le C_\delta\,\int_\Rd \phi(\rho,\ww)\,d\gamma\;.
\]
We can then pass to the limit as $\delta\down 0$, since $\rho_\delta\to\rho$ pointwise.

The last step is to approximate the functions $\rho$ and $\ww$, with $\delta\le \rho\le k$, $|\ww|\le k$, by smooth functions. We consider a family of smooth approximations $\tilde\rho_\eps$ and $\ww_\eps$ obtained by convolution with a smooth kernel. We finally set $m_\eps:=\int_\Rd \tilde\rho_\eps\,d\gamma$ and, in this framework, redefine $\rho_\eps:=m\,\tilde\rho_\eps/m_\eps$. Since $(\rho_\eps,\ww_\eps)$ converges to $(\rho,\ww)$ pointwise a.e.~in $\Rd$ and is uniformly bounded, we can pass to the limit as above when $\eps\down0$.
\end{proof}
\subsection{The action functional on measures}

Since we assumed that $h$ is concave and strictly positive for $\rho>0$, $h$ is an increasing map, so that $g$ is decreasing. We extend $h$ and $g$ to $[0,\infty)$ by continuity and we still denote by $\phi$ the lower semicontinuous envelope of $\phi$ in the closure $[0,\infty)\times \Rd$. If $h(0)>0$ then $g(0)<\infty$ and $\phi(0,\ww)=g(0)\,|\ww|^2$. When $h(0)=0$ we have $g(0)=\infty$ and
\[
\phi(0,\ww)=\begin{cases}
\infty&\mbox{if}\;\ww\neq 0\;,\\
0&\mbox{if}\;\ww=0\;.
\end{cases}
\]
We also introduce the \emph{recession functional}
\[
\label{eq:scheme:4}
\phi^\infty(\rho,\ww):=\sup_{\lambda >0}\frac 1\lambda\,\phi(\lambda\,\rho,\lambda\,\ww)=\lim_{\lambda \uparrow\infty}\frac 1\lambda\,\phi(\lambda\,\rho,\lambda\,\ww)\;,
\]
which is still a convex and lower semicontinuous function with values in $[0,\infty]$, and $1$-homogeneous. It is determined by the behaviour of $h(\rho)$ as $\rho\uparrow\infty$. If we set
\[
\label{eq:scheme:5}
h^\infty:=\lim_{\rho\uparrow\infty}\frac{h(\rho)}{\rho}=:\frac 1{g^\infty}\;,
\]
we have
\[
\label{eq:scheme:6}
\phi^\infty(\rho,\ww)=\begin{cases}
\infty&\quad\mbox{if}\;\ww\neq 0\\
0&\quad\mbox{if}\;\ww=0
\end{cases}
\quad\mbox{when}\quad h^\infty=0\;,
\]
and
\[
\label{eq:scheme:6bis}
\phi^\infty(\rho,\ww)=\begin{cases}
\frac{|\ww|^2}{h^\infty \rho}=g^\infty\frac{|\ww|^2}{\rho}&\quad\mbox{if}\;\rho\neq 0\\
\infty&\quad\mbox{if}\;\rho=0\;\mbox{and}\;\ww\neq 0
\end{cases}
\quad\mbox{when}\quad h^\infty>0\;.
\]

Let $\mu\in\FPM(\Rd)$ be a nonnegative Radon measure and let
$\nnu\in\mathcal M(\Rd;\Rd)$ be a vector Radon measure on $\Rd$. We write their Lebesgue decomposition with respect to the reference measure $\gamma$ as
\[
\label{eq:scheme:7}
\mu:=\rho\,\gamma+ \mu^\perp\,,\quad\nnu:=\ww\,\gamma+ \nnu^\perp\,.
\]
We can always introduce a nonnegative Radon measure $\sigma\in \FPM(\Rd)$ such that $\mu^\perp =\rho^\perp\sigma\ll\sigma$, $\nnu^\perp= \ww^\perp\sigma\ll\sigma$, e.g.~$\sigma:=\mu^\perp+|\nnu^\perp|$ and define the \emph{action functional}
\[
\label{eq:scheme:8}
\Phi(\mu,\nnu\,|\,\gamma):=\int_{\Rd}\phi(\rho,\ww)\,d\gamma+\int_{\Rd} \phi^\infty(\rho^\perp,\ww^\perp)\,d\sigma\;.
\]
Since $\phi^\infty$ is $1$-homogeneous, this definition is independent of $\sigma$. As we have done up to now, we shall simply write $\Phi(\mu,\nnu)=\Phi(\mu,\nnu\,|\,\gamma)$ when there is no ambiguity on the reference measure $\gamma$.

\smallskip\begin{remark}
If $h$ has a \emph{sublinear growth}, then $h^\infty=0$ and, as a consequence, if $\Phi(\mu,\nnu)<+\infty$, then we have
\[
\label{eq:scheme:16}
\nnu=\ww\cdot \gamma\ll\gamma\quad\mbox{and}\quad\Phi(\mu,\nnu)=\int_{\Rd} \phi(\rho,\ww)\,d\gamma\;,
\]
so $\Phi(\mu,\nnu)$ is independent of the singular part $\mu^\perp$. When $h$ has a \emph{linear growth}, i.e.~$h^\infty>0$, if $\Phi(\mu,\nnu)<+\infty$, then we have
\[
\nnu^\perp=\ww^\perp\cdot \mu^\perp\ll\mu^\perp\,.
\]
In both cases, one can choose $\sigma=\mu^\perp$, so that
\[
\nnu=\ww\cdot \gamma\ll\gamma+\ww^\perp\cdot \mu^\perp
\]
and if $g^\infty=1/h^\infty$ is finite, then we have
\[
\label{eq:scheme:17}
\Phi(\mu,\nnu\,|\,\gamma)=\int_{\Rd} \phi(\rho,\ww)\,d\gamma+g^\infty\int_{\Rd} |\ww^\perp|^2\,d\mu^\perp\,,
\]
while the last term simply drops if $h^\infty=0$.
\end{remark}

\smallskip\begin{lemma}[Lower semicontinuity, regular approximation of the action functional]\label{le:action_lsc}
The action functional is lower semicontinuous with respect to the weak convergence of measures, i.e.~if $(\gamma_n)$, $(\mu_n)$ and $(\nu_n)$ are sequences such that $\gamma_n\weakto\gamma$ weakly in $\FPM(\Rd)$, $\mu_n\weakto\mu$ weakly in $\FPM(\Rd)$ and $\nnu_n\weaksto\nnu$ in $\mathcal M(\Rd;\Rd)$ as $n\up\infty$, then
\[
\liminf_{n\uparrow \infty}\Phi(\mu_n,\nnu_n\,|\,\gamma_n)\ge\Phi(\mu,\nnu\,|\,\gamma)\;.
\]
Moreover, for every $\mu\in\FPM(\Rd)$ and $\nnu\in\mathcal M(\Rd;\Rd)$ such that $\Phi(\mu,\nnu)<\infty$, there exist sequences $(\mu_n)$ and $(\nu_n)$ for which
\[
\label{eq:scheme3:6}
\mu_n:=\rho_n\,\gamma\;\mbox{with}\;\rho_n\in C^0_b(\Rd)\;\mbox{and}\;\inf \rho_n>0\;,\quad\nnu_n:=\ww_n\,\gamma\;\mbox{with}\;\ww_n\in C^0_b(\Rd;\Rd)
\]
such that
\begin{equation}
\label{eq:scheme3:7}
\mu_n\weakto\mu\;\mbox{and}\;\nnu_n\weakto\nnu\;,\quad\lim_{n\uparrow\infty}\int_\Rd \phi(\rho_n,\ww_n)\,d\gamma=\Phi(\mu,\nnu\,|\,\gamma)\;.
\end{equation}
\end{lemma}\smallskip
\begin{proof}
The first statement is a well known fact about lower semicontinuity of convex integrals (see e.g.~\cite{ambrosio-Fusco-Pallara00}). Concerning the approximation property \eqref{eq:scheme3:7}, general relaxation results provide a family of approximations in $L^1_\gamma(\Rd)$. We can then apply Lemma \ref{le:smooth_approximations} and a standard diagonal argument.
\end{proof}

\subsection{The weighted Wasserstein distance}\label{subsec:WeightedWasserstein}

Denote by $\mathcal B(\Rd )$ the collection of all Borel subsets of $\Rd$, by $\FPM(\Rd)$ the collection of all finite positive Borel measures defined on $\Rd$ and by $\Probabilities\Rd$ the convex subset of all probability measures i.e.~all $\mu\in\FPM(\Rd)$ such that $\mu(\Rd)=1$. If $\mathcal M(\Rd;\Rd)$ is the set of the vector valued Borel measures $\nnu:\mathcal B(\Rd) \to\Rd$ with finite variation, i.e.~such that
\begin{multline*}
\label{eq:scheme3:16ter}
|\nnu|(B):=\sup\Big\{\sum_{j\le n} |\nnu(B_j)|\,:\\
B=\bigcup_{j\le n}B_j\,,\;B_j\in\mathcal B(\Rd)\;\mbox{pairwise disjoint}\,,\;n<\infty\Big\} <\infty
\end{multline*}
for any $B\in\mathcal B(\Rd)$, then $|\nnu|$ is in fact a finite positive measure in $\FPM(\Rd)$ and $\nnu$ admits the polar decomposition $\nnu=\ww\,|\nnu|$ where the Borel vector field $\ww$ belongs to $L^1_{|\snnu|}(\Rd;\Rd)$. We can also consider $\nnu$ as a vector $(\nnu^1,\nnu^2,\cdots,\nnu^d)$ of $d$ measures in $\mathcal M(\Rd;\R)$.

For any $T>0$, let $\mathcal{CE}(0,T;\Rd)$ be the set of time dependent measures $(\mu_t)_{t\in [0,T]}$, $(\nnu_t)_{t\in (0,T)}$ such that
\begin{enumerate}
\item $t\mapsto \mu_t$ is weakly $^*$ continuous in $\RPM(\Rd)$,\vspace*{2pt}
\item $(\nnu_t)_{t\in (0,T)}$ is a Borel family with $\int_0^T |\nnu_t|(B_R)\,d t<\infty$ for any $R>0$,\vspace*{6pt}
\item $(\mu,\nnu)$ is a distributional solution of
\[
\label{eq:continuity1}
\partial_t {\mu}_t+\nabla\cdot \nnu_t=0\quad\mbox{in}\;\R^d\times (0,\FinalT)\;.
\]
\end{enumerate}
As in \cite{MR2448650}, we define the weighted Wasserstein distance as follows.
\smallskip\begin{definition}
\label{def:main}
The $(h,\gamma)$-Wasserstein distance between $\mu_0$ and $\mu_1\in\RPM(\Rd)$ is defined by
\begin{multline}
\label{eq:10}
W_{h,\gamma}(\mu_0,\mu_1):=\inf\Big\{\Big[\textstyle{\int_0^1\Phi(\mu_t,\nnu_t\,|\,\gamma)\,dt}\Big]^{1/2}\,:\\
(\mu,\nnu)\in\mathcal{CE}(0,1;\Rd)\,,\;\mu_{t=0}=\mu_0\,,\;\mu_{t=1}=\mu_1\Big\}
\end{multline}
with $\Phi(\mu,\nnu\,|\,\gamma):=\Phi(\rho,\ww)+\Phi^\infty(\ww^\perp)$ if $\mu=\rho\,\gamma+\mu^\perp$ and $\nnu=\ww\,\gamma+\ww^\perp\,\mu^\perp$, $\Phi(\mu,\nnu\,|\,\gamma):=\infty$ otherwise, and $\Phi^\infty(\ww):=\lim_{\lambda \uparrow\infty}\lambda\,\phi(\lambda,\ww)$.

We denote by $\mathcal M_{h,\gamma}[\sigma]$ the set of all measures $\mu\in\RPM(\Rd)$ which are at finite $W_{h,\gamma}$-distance from $\sigma$.
\end{definition}\smallskip

Notice that in \cite{MR2448650} we were using the notation
$W_{\phi,\gamma}$ instead of $W_{h,\gamma}$. Whenever there is no
ambiguity on the choice of the measure $\gamma$, we shall simply write
$W_h$. The next result is taken from
\cite[Theorem 5.6 and Proposition 5.14]{MR2448650}
\smallskip\begin{theorem}[Lower semicontinuity]
\label{thm:2}
If $\phi$ satisfies \eqref {eq:2} and \eqref{eq:53}, the map
$(\mu_0,\mu_1)\mapsto W_{h,\gamma}(\mu_0,\mu_1)$ is lower
semicontinuous with respect to the weak $^*$~convergence in
$\RPM(\Rd)$. More generally, suppose that $\gamma^n\weaksto \gamma$ in
$\RPM(\Rd)$,
$h^n$ is monotonically decreasing w.r.t.~$n$ and pointwise
converging to $h$, and $\mu^n_0\weaksto\mu_0$, $\mu^n_1\weaksto\mu_1$ in $\RPM(\Rd)$ as $n\up\infty$. Then
\[
\label{eq:cap5:9}
\liminf_{n\up\infty}W_{h^n,\gamma^n}(\mu^n_0,\mu^n_1)\ge W_{h,\gamma}(\mu_0,\mu_1)\;.
\]
If moreover $\gamma_n\ge\gamma$ we have
\[
\label{eq:stability-h}
\lim_{n\to+\infty}W_{h_n,\gamma_n}(\mu_0,\mu_1) = W_{h,\gamma}(\mu_0,\mu_1)\;.
\]
\end{theorem}\smallskip
It is possible to reparametrize the path connecting $\mu_0$ to $\mu_1$ in the definition of $W_{h,\gamma}$ and establish that, for any $T>0$,
\[
W_{h,\gamma}(\sigma,\eta):=\inf\Big\{\sqrt T\,\Big[\textstyle{\int_0^T\Phi(\mu_t,\nnu_t\,|\,\gamma)\,dt}\Big]^{1/2}\,:(\mu,\nnu)\in{\mathcal CE}(0,T;\sigma\to\eta)\Big\}
\]
where ${\mathcal CE}(0,T;\sigma\to\eta)$ denotes the set of the paths $(\mu,\nnu)\in\mathcal{CE}(0,T;\Rd)$ such that $\mu_{t=0}=\sigma$ and $\mu_{t=T}=\eta$. By \cite[Theorem 5.4 and Corollary 5.18]{MR2448650}, we have the
\smallskip\begin{theorem}[Existence of geodesics]
\label{existence:geodesics}
Whenever the infimum in \eqref{eq:10} has a finite value, it is attained by a curve $(\mu,\nnu)\in\CE\phi(0,1;\Rd)$ such that
\[
\label{eq:cap3:37bis}
\Phi(\mu_t,\nnu_t\,|\,\gamma)=W_{h,\gamma}^2(\mu_0,\mu_1)\quad\forall\;t\in (0,1)\;\Leb 1\,\mbox{a.e.}
\]
In this case we have the equivalent characterization
\[
\label{eq:cap3:29}
W_{h,\gamma}(\sigma,\eta)=\min\Big\{{\textstyle\int_0^T \big[\Phi(\mu_t,\nnu_t\,|\,\gamma)\big]^{1/2}dt}\,:\,(\mu,\nnu)\in{\mathcal CE}(0,T;\sigma\to\eta)\Big\}\,.
\]
The curve $(\mu_t)_{t\in [0,1]}$ associated to a minimum for~\eqref{eq:10} is a constant speed mimimal geodesic:
\[
\label{eq:cap5:7}
W_{h,\gamma}(\mu_s,\mu_t)=|t-s|\,W_{h,\gamma}(\mu_0,\mu_1)\quad\forall\;s\,,\;t\in [0,1]\;.
\]
\end{theorem}\smallskip
We may notice that the characterization of $W_{h,\gamma}(\sigma,\eta)$ in terms of $\int_0^T \sqrt{\Phi(\mu_t,\nnu_t\,|\,\gamma)}\,dt$ allows to consider the case $T=+\infty$. By \cite[Chap. 1]{Ambrosio-Gigli-Savare05}
(also see \cite[p. 222]{MR2448650}), one knows that
\[
W_{h,\gamma}(\mu_0,\mu_T) \leq \int_0^T|\mu_t'|\,dt\quad\mbox{with}\;|\mu_t'|:=\lim_{h\to 0}\frac{W_{h,\gamma}(\mu_{t+h},\mu_t)}h
\]
for any absolutely continuous curve $t\mapsto\mu_t$ such that $\mu_{t=0}=\mu_0$ and $\mu_{t=T}=\mu_T$.

\medskip Now let us come back to the formal point of view of Section~\ref{Sec:Formal} and establish \eqref{eq:72} in this framework. Assume that $\rho_t$ is given by KFP and $\ww_t=\dD\rho_t$. The curve $\mu_t=\rho_t\,\gamma$ connects $\mu_0=\rho_0\,\gamma$ with $\mu_\infty=\gamma$ and, using
\[
\sqrt{P(\rho_t)}=\sqrt{\Phi(\rho_t,\ww_t\,|\,\gamma)}=|\mu_t'|\;,
\]
it follows that
\[
W_{h,\gamma}(\rho_0,\gamma)\le\int_0^{\infty}\sqrt{P(\rho_t)}\,dt=\int_0^{\infty}|\dot\mu_t|\,dt
\]
as already noted in Section~\ref{Sec:Formal} (equality case in \eqref{eq:69bis}). On the other hand, for any $(\mu,\nnu)\in{\mathcal CE}(0,T;\mu_0\to\mu_T)$, $T\in(0,\infty)$, we have $|\dot\mu_t|\le\sqrt{\Phi(\mu_t,\nnu_t)}$ and so
\[
\int_0^T|\dot\mu_t|\,dt\le\int_0^T\sqrt{\Phi(\mu_t,\nnu_t)}\,dt\;.
\]
By taking first the infimum $(\mu,\nnu)\in{\mathcal CE}(0,T;\mu_0\to\mu_T)$ and then the limit $T\to\infty$, we also find
\[
\int_0^{\infty}|\dot\mu_t|\,dt\le W_{h,\gamma}(\rho_0,\gamma)\;,
\]
thus proving the equality in the above inequality. This completes the proof of~\eqref{eq:72}.

\section{Entropy and entropy production}\label{subsec:entropy}

Let us consider now a function $\psi$ such that $\psi''(x)=g(x)$ for any $x>0$. Among all possible choices of $\psi$, we consider in particular the convex functions $\psi_a:[0,\infty)\to[0,\infty)$ depending on $a>0$ and characterized by the conditions
\[
\label{eq:scheme3:16}
\psi_a''(x)=g(x)\;,\quad \psi_a(a)=\psi_a'(a)=0\;,\quad\mbox{i.e.}\quad\psi_a(x)=\int_a^x (x-r)\,g(r)\,dr\;.
\]
Observe that $\psi_a\in C^2(0,\infty)$ has a strict minimum at $a>0$ and it satisfies the transformation rule
\[
\label{eq:scheme3:36}
\psi_a(x)=\psi(x)-\,\psi(a)-\,\psi'(a)\,(x-a)\quad\forall\;a>0\;,
\]
independently of the choice of $\psi$ (for a given function $g$). When $g(x)=1/x$ we obtain the logarithmic entropy density $\sfE(x):=x\log x$ and the family
\[
\begin{aligned}
\sfE_a(x):=\int_a^y (y-r)\,\frac1r\,dr=x\log x-a\log a -(1+\log a)\,(x-a)\;,
\end{aligned}
\label{eq:scheme3:17}
\]
which provides useful lower/upper bounds for $\psi$. In fact, $h$ being concave, if $h(0)=0$, then $h(x)\ge h(a)\,x$ if $0<x<a$, so that
\[
\label{eq:scheme3:16bis}
g(x)\le \frac {g(a)}x\quad\mbox{and}\quad\psi(x)\le g(a)\,\sfE_a(x)\quad\forall\;x\in(0,a]\;.
\]
On the other hand, when $x\ge a$, we have $h(x)\le h(a)\,x$, so that
\begin{equation}
\label{eq:scheme3:12}
g(x)\ge\frac {g(a)}x\quad\mbox{and}\quad\psi(x)\ge g(a)\,\sfE_a(x)\quad\forall\;x\in[a,+\infty)\;,
\end{equation}
thus showing that $\psi(x)$ has a superlinear growth as $x\uparrow\infty$.

We can therefore introduce the \emph{relative entropy functional}
\[
\label{eq:scheme3:37}
\Psi(\rho):=\int_\Rd\psi_a(\rho(x))\,d\gamma(x)=\int_\Rd\Big(\psi(\rho(x))-\,\psi(a)\Big)\,d\gamma\quad\mbox{with}\quad a=\int_\Rd\rho\,d\gamma\;.
\]
In the particular case $\psi=\sfE$, we set
\[
\label{eq:Cap2:10}
\mathcal H(\rho):=\int_\Rd\rho\log\rho\,d\gamma-a\log a\quad\quad\mbox{with}\quad a=\int_\Rd\rho\,d\gamma\;.
\]
Since $\psi$ is convex and superlinearly increasing, if $\sup_n\Psi(\rho_n)<\infty$, then there exists a subsequence weakly converging to $\rho$ in $L^1_\gamma(\Rd)$ and
\[
\label{eq:scheme3:39}
\liminf_{n\uparrow\infty}\Psi(\rho_n)\ge \Psi(\rho)\;.
\]
\smallskip\begin{remark}
\label{le:McCann}
If the function $\psi$ satisfies $\psi''=g$, $\psi(0)=0$ and if~\eqref{eq:2} holds, then $\psi$ also satisfies McCann's conditions, i.e.~the map $x\mapsto e^{x}\,\psi(e^{-x})$ is convex and non increasing on $(0,\infty)$ or, equivalently,
\[
\label{eq:cap3:5}
x\,\psi'-\,\psi\ge 0\quad\mbox{and}\quad x^2\,\psi''-x\,\psi'+\psi\ge0\quad\forall\;x>0\;.
\]

The convexity of $\psi$ indeed yields $x\,\psi'(x)-\,\psi(x)\ge -\,\psi(0)=0$. Consider the function $\vartheta(x):= x^2\,\psi''(x)-x\,\psi'(x)+\psi(x)$ and observe that $\lim_{x\down0}\vartheta(x)=0$, since $\psi''=1/h$ and $h$ is concave so that, in particular, $h(x)\ge c\,x$ near $x=0$, for some positive constant $c$. On the other hand, we have
\[
\vartheta'(x)=x^2\,g'(x)+x\,g(x)=x\,\frac d{dx}\left(\frac x{h(x)}\right)
\]
and the function $x\mapsto h(x)/x$ being positive, non increasing, we deduce that \hbox{$\vartheta'(x)\ge0$}, so that $\vartheta\ge0$.
\end{remark}

\medskip Let us introduce the Sobolev spaces
\[
\label{eq:Cap2:11}
W^{1,p}_\gamma(\Rd):=\Big\{\rho\in W^{1,p}_{\rm loc}(\Rd)\,:\,\int_\Rd\Big(|\rho|^p+|\dD\rho|^p\Big)\,d\gamma<\infty\Big\}\,.
\]
For $\rho\in W^{1,1}_\gamma(\Rd)$, $\rho\ge0$, we define the \emph{entropy production functional} as
\[
\label{eq:scheme3:21}
P(\rho):=\Phi(\rho,\dD\rho)\;\mbox{with domain}\;\mathcal D(P):=\Big\{\rho\in W^{1,1}_\gamma(\Rd)\,:\,\rho\ge0\;,\quad P(\rho)<\infty\Big\}\,.
\]
We also introduce the absolutely continuous functions
\[
\label{eq:Cap2:2}
f(r):=\int_0^r \sqrt {g(\xi)}\,d\xi\;,\quad L_\psi(r):=r\,\psi'(r)-\,\psi(r)\;,
\]
and observe that
\[
\label{eq:Cap2:15}
\frac{d}{dr}L_\psi(r)=r\,\psi''(r)=r\,g(r)=\frac r{h(r)}
\]
is bounded if and only if $h(r)$ has a linear growth as $r\up\infty$. In the case $h(r)=r$, $\psi=\sfE$, to the entropy functional $\mathcal H$ corresponds the \emph{entropy production functional}
\[
\label{eq:Cap2:18}
\mathcal I(\rho):=\int_\Rd\frac{|\dD\rho|^2}\rho\,d\gamma\;.
\]
\smallskip\begin{proposition}
\label{thm:production}
Let $\rho$ be nonnegative function in $L^1_\gamma(\Rd)$. Then $\rho\in W^{1,1}_\gamma(\Rd)$ and $P(\rho)<\infty$ if and only if $\dD f(\rho)\in L^2_\gamma(\Rd;\Rd)$ and in this case we have
\[
P(\rho)=\int_\Rd |\dD f(\rho)|^2\,d\gamma\;.
\]
If $\rho\in\mathcal D(P)$ and $h(r)\ge\sfh\,r$ for some constant $\sfh>0$, then $L_\psi(\rho)\in W^{1,1}_\gamma(\Rd)$,
\begin{equation}
\label{eq:Cap2:16}
\int_\Rd\frac{|\dD L_\psi(\rho)|^2}\rho\,d\gamma\le\sfh^{-1}P(\rho)\quad\mbox{and}\quad P(\rho)\le \sfh^{-1}\,\mathcal I(\rho)\;.
\end{equation}
Moreover, the functional $\rho\mapsto P(\rho)$ is lower semicontinuous with respect to the weak convergence in $L^1_\gamma(\Rd)$, i.e.~if a sequence $(\rho_n)_{n\in\N}$ weakly converges to some $\rho$ in $L^1_\gamma(\Rd)$ and $\sup_{n\in\N}P(\rho_n)<\infty$, then $\rho\in W^{1,1}_\gamma(\Rd)$ and
\begin{equation}
\label{eq:Cap2:4}
\liminf_{n\uparrow\infty}P(\rho_n)\ge P(\rho)\;.
\end{equation}
\end{proposition}\smallskip
\begin{proof}
Identity \eqref{eq:Cap2:4} and $\dD\,L_\psi(\rho)=\rho\,g(\rho)\,\dD\rho$ are straightforward if $\rho$ takes its values in a compact interval of $(0,\infty)$. The general case follows as in Lemma~\ref{le:smooth_approximations} by a standard truncation argument, while the lower semicontinuity is a consequence of convexity.
\end{proof}

\section{The KFP flow and its first variation}\label{Sec:Variational}

\subsection{Variational solutions to the KFP flow}\label{KFP}

As in Section~\ref{Sec:Formal}, let us introduce the differential operators
\[
\label{eq:scheme:27}
\begin{aligned}
\tdiv \vv:=&\,e^V\,\nabla\cdot (e^{-V}\vv)=\nabla\cdot\vv- \vv\cdot \dD V\;,\\
\tDelta \rho:=&\,\tdiv(\dD\rho)=\Delta \rho-\dD\rho\cdot \dD V\;,
\end{aligned}
\]
which, with respect to the measure $\gamma$, satisfy the following ``integration by parts formulae'' against test functions $\zeta\in C^\infty_c(\Rd)$:
\[
\label{eq:233}
\begin{aligned}
\int_{\Rd} \vv\cdot \dD\zeta\,d\gamma=-\int_{\Rd} \tdiv\vv\,\zeta\,d\gamma\quad\mbox{and}\quad\int_{\Rd} \dD v \cdot \dD\zeta\,d\gamma=-\int_{\Rd} \tDelta v\,\zeta\,d\gamma\;.
\end{aligned}
\]
We consider the Kolmogorov-Fokker-Planck equation
\begin{equation}
\label{eq:cap4:3}
\partial_t\rho_t-\tDelta\rho_t=0\quad\mbox{in}\;(0,\infty)\times\Rd\;.
\end{equation}
For simplicity, we will consider equations in the whole $\R^d$
(corresponding to the finiteness assumption on the potential $V$);
necessary adaptations when this is not the case are straightforward
and left to the reader.
We will also assume that
\begin{equation}
 \label{eq:1}
 \text{\emph{the potential $V$ is smooth with
bounded second derivatives\,.}}
\end{equation} 
Based on the integration by parts formula, the variational formulation of \eqref{eq:cap4:3} in the Hilbert space $L^2_\gamma(\Rd)$ relies on the symmetric, closed Dirichlet form
\[
\label{dirigamma}
a_\gamma(\rho,\eta):=\int_{\Rd}\langle\dD\rho,\dD\eta\rangle\,d\gamma\quad\forall\;\rho\,,\;\eta\in W^{1,2}_\gamma(\Rd)\;,
\]
where $W^{1,2}_\gamma(\Rd)$ is endowed with its natural norm $\|\rho\|_{W^{1,2}_\gamma(\Rd)}^2:=\|\rho\|_{L^2_\gamma(\Rd)}^2+\,a_\gamma(\rho,\rho)$. Using smooth approximations, it is not difficult to prove that $W^{1,2}_\gamma(\Rd)$ is dense in $L^2_\gamma(\Rd)$. The abstract theory of variational evolution equation and the log-concavity of the measure $\gamma$ yield the following result (see e.g.~\cite[Thm. 6.7]{Ambrosio-Savare06}).
\smallskip\begin{proposition}
\label{prop:Fkpre}
Assume that \eqref{eq:Hessian}--\eqref{eq:2} hold. For every $\rho_0\in L^2_\gamma(\Rd)$, the solution of \eqref{eq:cap4:3} has the following properties:
\begin{enumerate}
\item
There exists a unique $\rho_t=S_t\rho_0\in W^{1,2}_{\rm loc}\left(0,\infty;L^2_\gamma(\Rd)\right)$, $t>0$, such that
\begin{equation}
\label{eq:Gibbs:16}
\frac{d}{dt}\langle\rho_t,\eta\rangle_{L^2_\gamma(\Rd)}+a_\gamma(\rho_t,\eta)=0\quad\forall\;\eta\in W^{1,2}_\gamma(\Rd)\;,\quad \lim_{t\downarrow0}\rho_t=\rho_0\;\mbox{in}\;L^2_\gamma(\Rd)\;.
\end{equation}
If $\rho_{\rm min}\leq\rho_0\leq\rho_{\rm max}$, then $\rho_t$ satisfies the same uniform bounds. The semigroup $(S_t)_{t\geq 0}$ is an analytic Markov semigroup in $L^2_\gamma(\Rd)$ which can be extended by continuity to a contraction semigroup in $L^p_\gamma(\Rd)$ for every $p\in [1,\infty)$ and to a weakly $^*$ continuous semigroup in $L^\infty_\gamma(\Rd)$.
\item
For every $\rho$, $\sigma\in L^2_\gamma(\Rd)$, we have
\[
\label{eq:Cap2:12}
\int_\Rd\big(S_t\rho\big)\,\sigma\,d\gamma=\int_\Rd\rho\,\big(S_t\sigma\big)\,d\gamma\quad\forall\;t\ge0\;.
\]
\item
For every $t>0$, $S_t$ maps $L^\infty_\gamma(\Rd)$ into $C_b(\Rd)$ and ${\rm Lip}_b(\Rd)$ into itself, with the uniform bound
\[
\label{feller}
[S_t \rho]_{{\rm Lip}(\Rd)}\leq [\rho]_{{\rm Lip}(\Rd)}\quad\forall\;t\geq 0\;,\quad\forall\;\rho\in {\rm Lip}_b(\Rd)\;.
\]
\item If $\rho_0\ge0$, $\int_\Rd|x|^2\,\rho_0\,d\gamma<\infty$ and $\mathcal H(\rho_0)<\infty$, then the map $t\mapsto \mathcal H(\rho_t)$ is convex, $\rho_t\in W^{1,1}_\gamma(\Rd)$ for every time $t>0$, and
\[
\label{eq:Cap2:23}
\sup_{t\in [0,T]}\int_\Rd|x|^2\rho_t\,d\gamma<\infty\;,\quad\frac d{dt}\,\mathcal H(\rho_t)=-\mathcal I(\rho_t)\;,\quad\frac d{dt}\Big(e^{2\lambda t}\,\mathcal I(\rho_t)\Big)\le 0\;.
\]
\end{enumerate}
\end{proposition}\smallskip
Notice that the Assumption $\rho_0\in L^2_\gamma(\Rd)$ is not needed in Property 4, according to~\cite[Thm. 6.7]{Ambrosio-Savare06}.

\subsection{Measure valued solutions to the FP flow}\label{KFPmeasure}

We first recall some basic results on measure-valued solutions of the Fokker-Planck (FP) equation
\begin{equation}
\label{fk}
\partial_t\mu_t = \Delta\mu_t + \nabla \cdot (\dD V\mu_t)\quad (t,x)\in(0,+\infty)\times\Rd\,.
\end{equation}
Solutions of \eqref{fk} are understood in the sense of distributions, i.e.~for any $T>0$ and $\varphi\in C^\infty_c([0,T]\times\Rd)$, we have
\begin{equation}
\label{eq:Gibbs:3}
\int_{\Rd}\varphi_T\,d\mu_T=\int_{\Rd}\varphi\,d\mu_0+\int_0^T\int_{\Rd}\Big(\partial_t\varphi_t+\Delta\varphi_t-\dD V\cdot\dD\varphi_t \Big)\,d\mu_t\,dt\;.
\end{equation}
For any $\mu\in\FPM(\Rd)$, we denote by $\sfm_p(\mu)$, $p\in [1,\infty)$, the $p$-moment of $\mu$, i.e. $\sfm_p(\mu):=\int_\Rd |x|^p\,d\mu(x)$. By $\FPMMM m2(\Rd)$ we denote the space of probability measures on~$\Rd$ with finite second moment $\sfm_2$. The \emph{relative entropy} of $\mu$ with respect to $\gamma$ is defined as
\[
\mathcal H(\mu\,|\,\gamma):=\int_{\Rd}\rho\,\log\rho\,d\gamma\;\mbox{if}\;\mu\ll\gamma\;\mbox{and}\;\mu=\rho\,\gamma\;,\quad\mathcal H(\mu\,|\,\gamma):=+\infty\;\mbox{otherwise}\;.
\]
Given two probability measures $\mu$ and $\nu$ in $\Probabilities\Rd$, the classical \emph{Wasserstein distance} $W_2$ is defined as $W_2(\mu,\nu):= \textstyle\inf\{[\int_{\Rd\times \Rd}|y-x|^2\,d\,\Sigma ]^{1/2}\,:\,\Sigma\in\Gamma(\mu,\nu)\}$. Here $\Gamma(\mu,\nu)$ is the set of all \emph{couplings} between $\mu$ and $\nu$: it consists of all probability measures $\Sigma$ on $\Rd\times \Rd$ whose first and second marginals are respectively $\mu$ and $\nu$, i.e.~$\Sigma(B\times \Rd)=\mu(B)$ and $\Sigma(\Rd\times B)=\nu(B)$ for any $B\in\BorelSets{\Rd}$. Notice that the notation $W_2$ is not consistent with the one for weighted distances $W_h$; we shall however use it as it is classical.

\medskip For a proof of the next results see e.g.~\cite[Sect. 3]{Ambrosio-Savare-Zambotti07}.
\smallskip\begin{proposition}[Uniqueness and stability of the solutions of FP]
\label{fpunique}
Let $\mu_0\in\FPMMM m2(\Rd)$.
\begin{enumerate}
\item The FP equation \eqref{eq:Gibbs:3} has a unique solution $\mu_t=\mathcal \mathcal S_t\mu_0$ in the class of weakly continuous maps $t\mapsto\mu_t\in\FPM (\Rd)$ with $\sup_{t\in(0,T)} \sfm_2(\mu_t)<+\infty$.
\item The unique solution $\mu_t$ is continuous with respect to the Wasserstein distance $W_2$ and Lipschitz continuous in all compact intervals $[t_0,t_1]\subset (0,+\infty)$.
\item It is characterized by the family of variational inequalities
\[
\label{eq:cap4:1}
\frac 12\,\frac d{dt}W_2^2(\mu_t,\nu)+\frac \lambda2\,W_2^2(\mu_t,\nu)+\mathcal H(\mu_t\,|\,\gamma)\le \mathcal H(\nu\,|\,\gamma)\quad\forall\;\nu\in\FPMMM m2(\Rd)\;.
\]
\item In addition, it is stable: $\mu^n_0\to\mu_0$ in $\FPMMM m2(\Rd)$ implies that $\mu^n_t\to\mu_t$ in $\FPMMM m2(\Rd)$ for all $t\geq 0$.
\end{enumerate}
\end{proposition}\smallskip

Notice that the measure $\gamma$ provides a \emph{stationary} solution of \eqref{fk}. All solutions $\mu_t$ weakly converge to $\gamma$ as $t\to +\infty$. Finally, $\mu_t$ is absolutely continuous with respect to $\gamma$ for any $t>0$, with density $\rho_t$, and $\rho_t$ is a solution of the KFP flow.

\subsection{Variational solutions to the modified KFP equation}

We consider the first variation of the KFP flow, i.e.~the \emph{modified Kolmogorov-Fokker-Planck equation}
\begin{equation}
\label{eq:34}
\partial_t \ww_t-\tDelta\ww_t+D^2V\,\ww_t=0\quad\mbox{in}\;(0,\infty)\times\Rd\;,\quad\lim_{t\downarrow0}\ww_t=\ww_0\quad\mbox{in}\;L^2_\gamma(\Rd;\Rd)
\end{equation}
for the vector field $\ww:(0,\infty)\times\Rd\to\Rd$.
In the Hilbert space $\wW:= W^{1,2}_\gamma(\Rd;\Rd)$, we consider the continuous (recall \eqref{eq:1}) bilinear form
\[
\label{eq:31}
\aa_\gamma(\vv,\ww):=\int_{\Rd} \Big(\dD\vv :\dD\ww+\dD ^2V\,\vv\cdot \ww\Big)\,d\gamma\;.
\]
We look for solutions $\ww\in W^{1,2}_{\rm loc} ((0,\infty);L^2_\gamma(\Rd;\Rd))\cap L^2_{\rm loc}([0,\infty);\wW)$ solving the variational formulation
\begin{equation}
\label{eq:32}
\frac d{dt}\int_{\Rd} \ww_t\cdot \zzeta\,d\gamma+\aa_\gamma(\ww_t,\zzeta)=0\quad\forall\;\zzeta\in\wW.
\end{equation}
Observe that vector fields in $C^1_c(\Rd;\Rd)$ belong to
$\wW$. Actually the space of smooth compactly supported functions
$C^\infty_c(\Rd;\Rd)$ is dense in $\wW$, and $\wW$ itself is dense in
$L^2_\gamma(\Rd;\Rd)$.
Notice moreover that if $\zeta:\R\to [0,\infty)$ is a smooth convex
function with bounded second order derivatives and
$\zeta(0)=0$, and $\zz(\ww):=\frac{\zeta'(|\ww|)}{|\ww|}\ww $ (with
$\zz(0)=0$), an easy calculation shows that solutions of \eqref{eq:32} satisfy
\begin{displaymath}
 -\frac {d}{dt}\int_\Rd \zeta(|\ww_t|)\,\d
 \gamma=\aa_\gamma(\ww_t,\zz(\ww_t))\ge0\quad\text{a.e.\ in } (0,\infty).
\end{displaymath}
With these observations in hand, we can apply the variational theory
of evolution equations
and a simple regularization argument to prove the next result. 
\smallskip\begin{proposition}
\label{prop:modified_KFP}
For every $\ww_0\in L^2_\gamma(\Rd;\Rd)$, there exists a unique solution $\ww=\rR\ww_0$ of \eqref{eq:32} in $W^{1,2}_{\rm loc} ((0,\infty);L^2_\gamma(\Rd;\Rd)) \cap L^2_{\rm loc}([0,\infty);\wW)$ with $\lim_{t\down0}\ww_t=\ww_0$ in $L^2_\gamma(\Rd;\Rd)$. The semigroup~$\rR$ is symmetric
\[
\label{eq:Cap2:13}
\int_\Rd\rR_t \ww\cdot \zz\,d\gamma=\int_\Rd\ww\cdot \rR_t\zz\,d\gamma\quad\forall\;\ww\,,\;\zz\in L^2_\gamma(\Rd;\Rd)\;,\quad\forall\;t>0\;,
\]
and satisfies
\[
 \label{eq:3}
 \int_\Rd \zeta\big(|\rR_t \ww_0|\big)\,\d\gamma\le \int_\Rd 
 \zeta\big(|\ww_0|\big)\quad
 \text{for every }\ww_0\in L^2_\gamma(\Rd;\Rd)
\]
and every convex function $\zeta:\R\to[0,\infty)$ with $\zeta(0)=0$.
In particular $\rR$ can be extended by density to a contraction
semigroup in $L^p_\gamma(\Rd;\Rd)$, $p\in [1,\infty]$.
\end{proposition}\smallskip

The link between \eqref{eq:cap4:3} and \eqref{eq:34} is enlightened by the next result.
\smallskip\begin{theorem}
\label{thm:density_gradient_link}
If $\rho_t$ is a variational solution of the KFP equation \eqref{eq:cap4:3} with initial datum $\rho_0\in W^{1,2}_\gamma(\Rd)$, then $\ww_t:=\dD\rho_t$ belongs to $C^0([0,\infty);L^2_\gamma(\Rd))$ and it is the solution of the modified KFP equation \eqref{eq:34} with initial datum $\ww_0:=\dD\rho_0$. In particular we have
\[
\label{eq:Cap2:14}
\int_\Rd \dD S_t\rho \cdot \ww\,d\gamma=\int_\Rd \dD\rho\cdot \rR_t\ww\,d\gamma\quad\forall\;\rho\in W^{1,2}_\gamma(\Rd)\;,\quad\forall\;\ww\in L^2_\gamma(\Rd;\Rd)\;.
\]
The same result holds if $\rho_0$ belongs to $W^{1,1}_\gamma(\Rd)$.
\end{theorem}\smallskip
\begin{proof}
Since $\mathcal D(\tDelta)$ is dense in $W^{1,2}_\gamma(\Rd)$, we can
assume that $\rho_0\in\mathcal D(\tDelta)$. Then the regularity result
of Proposition~\ref{prop:Fkpre} shows that $\rho_t\in\mathcal
D(\tDelta)$ for every $t\ge 0$. Setting $\ww_t:=\dD\rho_t$,
we know (see e.g.\ the argument in the proof of \cite[Lemma
5.2]{Gianazza-Savare-Toscani09})
that $\aa_\gamma(\ww_t,\ww_t)\le \|\tDelta \rho_0\|_{L^2_\gamma(\Rd)}<+ \infty$. For a fixed $\zzeta\in C^\infty_c(\Rd;\Rd)$, we can then evaluate
\begin{equation}
\label{eq:39}
\frac d{dt}\int_{\Rd} \ww_t\cdot \zzeta\,d\gamma=\frac d{dt}\int_{\Rd} \dD\rho_t\cdot\zzeta\,d\gamma=-\frac d{dt}\int_{\Rd} \rho_t\,\tdiv \zzeta\,d\gamma=\int_{\Rd} \dD\rho_t\cdot \dD (\tdiv \zzeta)\,d\gamma\;.
\end{equation}
With the notations $\partial_i=\partial/\partial x_i$ and $\partial_{ij}=\partial^2/\partial x_i\partial x_j$ for $i$, $j=1$, $2$\ldots $d$, let us observe that
\[
\big(\dD\,\tdiv \zzeta\big)_j=\sum_i \dxj\big(\dxi\zeta_i-\zeta_i\,\dxi V\big)=\sum_i \dxij\zeta_i -\dxj\zeta_i\,\dxi V-\zeta_i\,\dxij V
\]
and
\[
\dD\rho_t\cdot \dD (\tdiv \zzeta)=\sum_{i,j} \dxj \rho_t\,\dxij\zeta_i -\dxj\rho_t\,\dxj\zeta_i\,\dxi V-\dxj\rho_t\,\zeta_i\,\dxij V\;.
\]
Inserting this expression in \eqref{eq:39} and integrating by parts the first term we get
\begin{align*}
&\hspace*{-1cm}\int_{\Rd} \dD\rho_t\cdot \dD (\tdiv \zzeta)\,d\gamma
\\&=\sum_{i,j}\int_{\Rd}\Big(\dxj \rho_t\,\dxij\zeta_i\,d\gamma -\sum_{i,j}\int_{\Rd} \Big(\dxj\rho_t\,\dxj\zeta_i\,\dxi V+\dxj\rho_t\,\zeta_i\,\dxij V\Big)\,d\gamma
\\&=\sum_{i,j}\int_{\Rd}\Big( -\dxij\rho_t\,\dxj\zeta_i +\dxi V\,\dxj \rho_t\,\dxj\zeta_i\Big)\,d\gamma\\
&\hspace*{4cm}-\sum_{i,j}\int_{\Rd} \Big(\dxj\rho_t\,\dxj\zeta_i\,\dxi V+\dxj\rho_t\,\zeta_i\,\dxij V\Big)\,d\gamma
\\&=-\sum_{i,j}\int_{\Rd}\Big( -\dxij\rho_t\,\dxj\zeta_i\,\dxj\rho_t\,\zeta_i\,\dxij V\Big)\,d\gamma
\\&=-\int_{\Rd} \Big(\dD\ww_t:\dD\zzeta+D^2 \ww_t\cdot\zzeta\Big)\,d\gamma=-\,\aa_\gamma(\ww_t,\zzeta)\;.
\end{align*}
Combined with \eqref{eq:39}, this shows that $\ww_t:=\dD\rho_t$ satisfies
the variational formulation of \eqref{eq:34}. The case of $\rho_0\in
W^{1,1}_\gamma(\Rd)$ follows by a standard approximation
procedure,
the fact that $\dD S_t \rho_0=\rR_t \dD\rho_0$, and
the $L^1_\gamma$-contraction property of $\rR$.
\end{proof}

\subsection{Measure valued solutions to the modified KFP equation}

Exactly like the (K)FP equation, the modified system can be extended to vector-valued measures initial data. To $\ww_t$, we associate the vector valued measures $\nnu_t:=\ww_t\,\gamma\in\mathcal M(\Rd;\Rd)$ which satisfy the system
\[
\label{eq:cap4:MVMFP}
\partial_t\nnu_t=\Delta\nnu_t+\nabla\cdot(\dD V\otimes\nnu_t)-D^2V\,\nnu_t\;,
\]
in the weak sense, i.e.
\begin{equation}
\label{eq:cap4:4}
\frac d{dt}\int_\Rd \zzeta\cdot d\nnu_t=\int_\Rd \Big(\Delta \zzeta- \dD\zzeta\otimes\dD V-D^2V\,\zzeta\Big)\cdot\,d\nnu_t\quad\forall\;\zzeta\in C^2_c(\Rd)\;.
\end{equation}
The semigroup can be extended to initial data which are vector valued measures with finite total variation using equi-integrability and moment estimates taken from \cite{MR2448650}.
\smallskip\begin{proposition}[Equi-integrability and moment estimates]
Let $\zeta$ be a nonnegative Borel function such that $\mu(\zeta^2)=\int_\Rd \zeta^2\,d\mu$ and $\gamma(\zeta^2)= \int_\Rd\zeta^2\,d\gamma$ are finite. If $\Phi(\mu,\nnu)<\infty$, we have
\[
\label{eq:Cap3:1}
\Big(\int_\Rd\zeta\,d\,|\nnu|\Big)^2 \le \Phi(\mu,\nnu)\,\gamma(\zeta^2)\,h\left(\mu(\zeta^2)/\gamma(\zeta^2)\right)\;.
\]
\end{proposition}\smallskip
In particular, for every Borel set $A\in\BorelSets\Rd$ we have
\begin{equation}
\label{eq:scheme3:5}
\Big(|\nnu|(A)\Big)^2\le\Phi(\mu,\nnu)\,\gamma(A)\,h\big(\mu(A)/\gamma(A)\big)
\end{equation}
which in particular yields ($\gamma(\Rd)=1$)
\[
\label{eq:scheme3:5bis}
\Big(|\nnu|(\Rd)\Big)^2\le\Phi(\mu,\nnu)\,h\big(\mu(\Rd)\big)\,.
\]
If moreover $\sfm_2(\mu)<\infty$, we can bound the first moment of $|\nnu|$ by
\begin{equation}
\label{eq:cap3:1}
\sfm_1(|\nnu|)=\int_\Rd |x|\,d\,|\nnu|\le\Big(\Phi(\mu,\nnu)\,\sfm_2(\gamma)\,h\big(\sfm_2(\mu)/\sfm_2(\gamma)\big)\Big)^{1/2}\,.
\end{equation}
\smallskip\begin{theorem}\label{thm:vv-KFP}
For every $\nnu_0\in \mathcal M(\Rd;\Rd)$ with $\sfm_1(|\nnu_0|)<+\infty$, there exists a unique solution $\nnu_t=\mrR_t\nnu_0$ in the class of weakly continuous maps $t\mapsto\nnu_t\in \mathcal M(\Rd;\Rd)$ with $\sup_{t\in [0,T]}\sfm_1(|\nnu_t|)<+\infty$, for every final time $T>0$. When $\nnu_0=\ww_0\,\gamma$ then $\mrR_t\nnu_0=\rR_t\,\ww_0\,\gamma$. The map $\nnu_0\mapsto \mrR\nnu_0$ is stable in the following sense: if
\[
\label{eq:cap4:5}
\nnu^n_0\weaksto\nnu_0\quad\text{weakly$^*$ in $\mathcal M(\Rd;\Rd)$ with }\sup_n\sfm_1(|\nnu^n_0|)<+\infty
\]
then $\mrR_t\nnu^n_0\weaksto \mrR_t\nnu_0$ in $\mathcal M(\Rd;\Rd)$.
\end{theorem}\smallskip

\begin{proof}
We divide the proof in three steps.

\smallskip\noindent\textbf{Step 1.} Let us first associate to $\nnu=\ww\,\gamma\in \mathcal M(\Rd;\Rd)$ the probability measure
\begin{equation}
\label{eq:cap4:8}
\upsilon:=\frac1M\,\frac 1{\sqrt{1+|x|^2}}\,|\nnu|=\frac1M\,\frac{|\ww|}{\sqrt{1+|x|^2}}\,\gamma\in \Probabilities\Rd\,,
\end{equation}
where the constant $M$ is a renormalization factor such that $\upsilon(\Rd)=1$. Observe that if $\sfm_1(|\nnu|)=\int_\Rd |x|\,|\ww(x)|\,d\gamma$ is finite, then $\upsilon\in \ProbabilitiesTwo{\Rd}$ and
\begin{equation}
\label{eq:cap4:9}
\sfm_2(\upsilon)\le\frac1M\,\sfm_1(|\nnu|)\;.
\end{equation}
We also choose the action density to be $\phi_2(\rho,\ww):=|\ww|^2/\rho$ corresponding to $h(\rho)=\rho$, and observe that the corresponding functional writes
\begin{equation}
\label{eq:cap4:11}
\Phi_2(\upsilon,\ww)=M\int_\Rd \sqrt{1+|x|^2}\,|\ww(x)|\,d\gamma\le M\Big(|\nnu|(\Rd)+\sfm_1(|\nnu|)\Big)\;.
\end{equation}
\smallskip\begin{proposition}\label{prop:moment_estimate} Let us suppose that $\ww\in L^1_\gamma(\Rd;\Rd)$ with $\sfm_1(|\nnu|)<+\infty$ and set $\nnu:=\ww\,\gamma$, $\upsilon$ as in \eqref{eq:cap4:8}, $\ww_t=\rR_t\ww$, $\nnu_t=\ww_t\,\gamma$, $\upsilon_t=\mathcal \mathcal S_t\upsilon$. Then
 \begin{equation}
\label{eq:cap4:12}
|\nnu_t|(\Rd)\le|\nnu|(\Rd)\;,\quad\sfm_1(|\nnu_t|)\le
|\nnu|(\Rd)+2\,\sfm_1(|\nnu|)+4\,M\,\sfm_2(\gamma),
\end{equation}
and for any $t>0$, we have
\begin{equation}
\label{eq:cap4:14}
\Big(|\nnu_t|(A)\Big)^2\le M\Big(|\nnu|(\Rd)+\sfm_1(|\nnu|)\Big)\,\upsilon_t(A)\quad\forall\;A\in \BorelSets{\Rd}\;.
\end{equation}
\end{proposition}\smallskip
\begin{proof}
The first inequality of \eqref{eq:cap4:12} follows
by the $L^1_\gamma$-contraction property of $\rR$.
Since the FP flow contracts the Wasserstein distance by
Proposition~\ref{fpunique} and since $\gamma$ is a stationary
solution,
the triangle inequality for the Wasserstein distance and the fact
that $\sqrt{\sfm_2(\mu)}=W_2(\mu,\delta_0)$ yield 
\begin{multline*}
\sqrt{\sfm_2(\upsilon_t)}\le W_2(\upsilon_t,\gamma)+\sqrt{\sfm_2(\gamma)}
\le W_2(\upsilon,\gamma)+\sqrt{\sfm_2(\gamma)}\le\sqrt{\sfm_2(\upsilon)}+2\sqrt{\sfm_2(\gamma)}\;.
\end{multline*}
On the other hand, \eqref{eq:cap3:1} yields
\[
\sfm_1(|\nnu_t|)\le\sqrt{\sfm_2(\upsilon_t)}\,\sqrt{\Phi_2(\upsilon_t,\nnu_t)}\le\Big(\sqrt{\sfm_2(\upsilon)}+2\sqrt{\sfm_2(\gamma)}\Big)\,\sqrt{\Phi_2(\upsilon,\nnu)}\;.
\]
Here we used the fact that $\Phi_2(\upsilon_t,\nnu_t)\le\Phi_2(\upsilon,\nnu)$. This will appear later as a consequence of Theorem~\ref {thm:action_decay}, and is independent of the present result. Combined with \eqref{eq:cap4:9} and \eqref{eq:cap4:11}, this proves the estimate on $\sfm_1(|\nnu_t|)$.

Applying \eqref{eq:scheme3:5} and \eqref{eq:cap4:11}, we get \eqref{eq:cap4:14}.
\end{proof}

\smallskip\noindent\textbf{Step 2: existence.} Let us approximate a given $\nnu_0\in \mathcal M(\Rd;\Rd)$ with $\sfm_1(|\nnu_0|)<+\infty$ by a sequence $\nnu^k=\ww^k\,\gamma\weaksto \nnu_0$ as $k\to\infty$ in $\mathcal M(\Rd;\Rd)$ with $\ww_k\in L^2_\gamma(\Rd;\Rd)$ and $\sfm_1(|\nnu_k|)\to \sfm_1(|\nnu_0|)$. We set $\nnu^k_t:=\ww^k_t \,\gamma$ with $\ww^k_t=\rR_t \ww^k$ so that $\nnu^k_t$ solves \eqref{eq:cap4:4}. Thanks to Proposition~\ref{prop:moment_estimate}, we know that the first order moment of $\nnu^k_t$ are uniformly bounded. This is sufficient to pass to the limit (up to extraction of a suitable subsequence) in \eqref{eq:cap4:4} and to find a solution $\nnu_t$ which is weakly$^*$ continuous in $\mathcal M(\Rd;\Rd)$ and satisfies the initial condition in the sense that $\nu_t\weaksto \nu_0$ in $\mathcal M(\Rd;\Rd)$ as $t\downarrow 0$.

\smallskip\noindent\textbf{Step 3: uniqueness and stability.} It follows by a standard duality argument, like in the case of Equation \eqref{eq:34}. If $\nnu^1_t$ and $\nnu^2_t$ are two weakly continuous solutions of \eqref{eq:cap4:4}, their difference $\ssigma_t:=\nnu^1_t-\nnu^2_t$ solves
\begin{equation}
\label{eq:Gibbs:6}
\int_{\Rd}\zzeta_T\cdot d\ssigma_T=\int_0^T\int_{\Rd}\Big(\partial_t\zzeta_t+\Delta\zzeta_t - \dD\zzeta_t\otimes\dD V-D^2V \,\zzeta_t\Big)\cdot d\ssigma_t\,dt
\end{equation}
for every $T>0$ and $\zzeta\in C^\infty_{\rm c}([0,T]\times \Rd;\Rd)$. By a mollification technique, it is not difficult to check that \eqref{eq:Gibbs:6} also holds for every function $\zzeta\in C([0,T]\times\Rd;\Rd)$ with $\partial_t\varphi$, $\dD\zzeta$ and $D^2\zzeta$ continuous and bounded in $[0,T]\times\Rd$.

Next, we introduce a family of smooth convex potentials $V_n$ with bounded derivatives of arbitrary orders, which satisfies a uniform Lipschitz condition
\[
|\dD V_n(x)-\dD V_n(y)|\le L\,|x-y|\quad \forall\;x,\,y\in \Rd\,,
\]
for some positive constant $L$ which is independent of $n$ and such that
\begin{displaymath}
V_n\to V\,,\quad\dD V_n\to \dD V\,,\quad D^2V_n\to D^2V\quad\text{pointwise as }n\to\infty\;.
\end{displaymath} 
For a given $\eeta\in C^\infty_c(\Rd;\Rd)$, we consider the solution $\zzeta_t$ of the time reversed (adjoint) parabolic equation
\[
\label{eq:Gibbs:8}
\partial_t\zzeta_t+\Delta\zzeta_t -\dD\zzeta_t\cdot\dD V-D^2V\,\dD\zzeta_t=0\quad\text{in $(0,T)\times\Rd$}\;,\quad\zzeta_T=\eeta\;.
\]
Using a maximum principle that can be found in \cite{john} and the
fact that the first and second order spatial derivatives of $\zzeta$
solve an analogous equation, standard parabolic regularity theory
shows that $\zzeta$ is sufficiently regular to be used as a test
function in \eqref{eq:Gibbs:6}
and satisfies the uniform bound (observe that the second and third derivatives of $V_n$ are still uniformly bounded)
\[
\label{eq:Cap3:7}
\sup_{t,x}|\zzeta_n|+|\dD\zzeta_n|\le C<+\infty\;.
\]
This leads to
\[
\label{eq:Cap3:6}
\Big|\int_\Rd\eeta \cdot\, d\ssigma_T\Big|\le C\int_0^T \int_\Rd \Big(|\dD V-\dD V_n|+|D^2V-D^2V_n|\Big)\,d\,|\ssigma_t|\,dt\;.
\]
Since the first order moment of $|\ssigma_t|$ is uniformly bounded, we can pass to the limit as $n\to\infty$ obtaining $\int_\Rd \eeta\cdot d\sigma_T=0$. As $\eeta$ is arbitrary, we conclude that $\nnu_T^1=\nnu_T^2$. The stability is then a simple consequence of uniqueness.
\end{proof}

\section{Action decay along the KFP flow and consequences}

We can prove now our main estimate, which is a refined version of Theorem~\ref{Thm:MainEstimate}, under the assumption that
\begin{equation}
\label{Eqn:RefinedConvexity}
\begin{array}{c}\mbox{\emph{The function $h$ is concave and, for some $\beta\in[0,1)$,}}\\[6pt]
\mbox{\emph{$(1-\beta)\,h\,h''+2\,\beta\,(h')^2\le0$ holds in the sense of distributions.}}\end {array}
\end{equation}
Notice that \eqref{Eqn:RefinedConvexity} is equivalent to \eqref{eq:Cap2:31} with $\beta:=(1-\alpha)/(1+\alpha)$.
\smallskip\begin{theorem}
\label{thm:action_decay}
Assume that \eqref{eq:Hessian} and \eqref{eq:PartitionFunctionFinite} hold, and let $(\rho,\ww)\in L^1_\gamma(\Rd,\R^+)\times L^1_\gamma(\Rd,\Rd)$ be such that $\Phi(\rho,\ww)<\infty$. If \eqref{Eqn:RefinedConvexity} is satisfied, then
\[
\label{eq:100}
\Phi(S_t\rho,\rR_t\ww)+2\,\beta\sum_{i=1}^d\int_0^t\Phi(S_s\rho\;,\partial_i\rR_s\ww)\,e^{2\lambda(s-t)}\,ds \leq e^{-2\lambda t}\,\Phi(\rho,\ww)\quad\forall\;t\ge 0\;.
\]
\end{theorem}\smallskip
\begin{proof}
We first prove the result with the additional assumptions that $0<\rho_{\rm min}\le \rho\le \rho_{\rm max}$ and $|\ww|\le w_{\rm max}$ $\gamma$ a.e. in $\Rd$. Assume that $h$ is of class $C^2(0,\infty)$. It follows that $\rho_t=S_t\rho$ and $\ww_t=\rR_t\ww$ satisfy the same bounds and, for all $t>0$, $\rho_t$, $\partial_t\rho_t\in W^{1,2}_\gamma(\Rd)$ and $\ww_t$, $\partial_t\ww_t\in W^{1,2}_\gamma (\Rd;\Rd)$. The function $\phi$ is of class $C^2$ in the strip
\[
\label{eq:cap4:7}
Q:=[\rho_{\rm min},\rho_{\rm max}]\times \{\ww\in\Rd\,:\,|\ww|\le w_{\rm max}\}
\]
and its differential $D\phi(\rho,\ww)$ can be decomposed as
\[\!
\label{eq:cap4:6}
D_\rho\phi(\rho,\ww)=g'(\rho)\,|\ww|^2\;,\quad D_\sww \phi(\rho,\ww)=2\,g(\rho)\,\ww\;.
\]
Since $g(\rho)$ and $g'(\rho)$ are bounded, the differential is also in $L^2_\gamma(\Rd;\R^{d+1})$. As a consequence, the time derivative of $t\mapsto\Phi(\rho_t,\ww_t)$ exists and
\[
\frac{d}{dt}\,\Phi(\rho_t,\ww_t) = \int_{\Rd}\Big( g'(\rho_t)\,|\ww_t|^2\partial_t\rho_t+2\,g(\rho_t)\,\ww_t\cdot\partial_t\ww_t\Big)\,d\gamma\;.
\]
In order to apply \eqref{eq:Gibbs:16} and \eqref{eq:32} we have to verify that all components of $D\phi(\rho_t,\kern -1pt\ww_t)$ are in $W^{1,2}_\gamma(\Rd)$. We have already seen that they are in $L^2(\Rd)$. Let us compute their $x$-derivative:
\begin{align*}
\dD\left(D_\rho\phi(\rho_t,\ww_t)\right) &= g''(\rho_t)|\,\ww_t|^2\,\dD\rho_t+2\,g'(\rho_t)\,\ww_t\cdot \dD\ww_t\;,\\
\dD\left(D_{\sww^i}\phi(\rho_t,\ww_t)\right) &=2\,g'(\rho_t)\,\ww_t^i\,\dD\rho_t+2\,g(\rho_t)\,\dD\ww_t^i\quad\mbox{for any}\;i=1\,,\;2\,,\ldots\,d\;.
\end{align*}
The above functions are in $L^2_\gamma(\Rd)$, since $g(\rho_t)$, $g'(\rho_t)$, $g''(\rho_t)$ and $\ww_t$ are bounded, so we get
\begin{multline*}
\frac{d}{dt}\,\Phi(\rho_t,\ww_t) = -\sum_{i=1}^d\int_{\Rd}\Scalar{D^2\phi(\rho_t,\ww_t)(\partial_i\rho_t,\partial_i\ww_t)}{(\partial_i\rho_t,\partial_i\ww_t)}\;d\gamma\\
-2\int_{\Rd}g(\rho_t)\,D^2V\,\ww_t\cdot\ww_t\;d\gamma\;.
\end{multline*}
Recalling \eqref{eq:scheme:14} and the convexity assumption on $V$, we find
\[
\label{eq:101}
\frac{d}{dt}\,\Phi(\rho_t,\ww_t) \leq -2\,\beta\sum_{i=1}^d\Phi(\rho_t,\partial_i\ww_t)-2\lambda\,\Phi(\rho_t,\ww_t)\;.
\]
It follows from Gronwall's lemma that for all $s\in(0,t)$,
\[
\label{eq:102}
e^{2\lambda t}\,\Phi(\rho_t,\ww_t) +2\,\beta\sum_{i=1}^d\int_s^t\Phi(\rho_r,\partial_i\ww_r)\,e^{2\lambda r}\,dr\leq e^{2\lambda s}\,\Phi(\rho_s,\ww_s)\;.
\]
The result follows by passing to the limit as $s\downarrow 0$ and recalling that $\phi$ is continuous and bounded on $Q$ and $\rho_s$, $\ww_s$ converge to $\rho$, $\ww$ as $s\down 0$ in $L^2_\gamma(\Rd)$ and $L^2_\gamma(\Rd;\Rd)$ respectively. The general result for an arbitrary concave function $h$ easily follows by approximating $h$ by a decreasing family of smooth concave functions in the interval $[\rho_{\rm min},\rho_{\rm max}]$. Finally, the general case $\rho\in L^1_\gamma(\Rd)$, $\ww\in L^1_\gamma(\Rd;\Rd)$, without upper and lower bounds, follows by approximation, using Lemmas \ref{le:smooth_approximations} and \ref{le:action_lsc}.
\end{proof}

\medskip We can extend the results of Theorem~\ref{thm:action_decay} to measure valued initial data.
\smallskip\begin{corollary}
\label{cor:action_decay2}
Assume that \eqref{eq:Hessian}--\eqref{eq:2} hold. Let $\mu\in \ProbabilitiesTwo\Rd$ and $\nnu\in \mathcal M(\Rd;\Rd)$ with $\sfm_1(|\nnu|)<+\infty$. Then for every $t>0$ we have $\mu_t=\mathcal S_t\mu=\rho_t\,\gamma$, $\nnu_t=\mrR_t\nnu=\ww_t\,\gamma$ with $\rho_t\in W^{1,1}_\gamma(\Rd)$, $\ww_t\in W^{1,1}_\gamma(\Rd;\Rd)$ if $\beta>0$, and
\[
\label{eq:Cap3:2}
\Phi(\rho_t,\ww_t)+2\beta\sum_{i=1}^{n}\int_0^t\Phi(\rho_s,\partial_i\ww_s)\,e^{2\lambda(s-t)}\,ds \leq e^{-2\lambda t}\,\Phi(\mu,\nnu\,|\,\gamma)\quad\forall\;t>0\;.
\]
\end{corollary}\smallskip
\begin{proof}
This follows directly from the measure formulation of the KFP flow (Proposition \ref{fpunique} and Theorem \ref{thm:vv-KFP}).
\end{proof}

\medskip Let us now consider the entropy functional $\Psi(\rho):=\int_\Rd \psi(\rho)\,d\gamma$, for a function $\psi$ as in Section~\ref{subsec:entropy}.
\smallskip\begin{theorem}
\label{thm:Ent-Prod}
Let $\rho\in\mathcal D(\Psi)$ with $\int_\Rd|x|^2\,\rho\,d\gamma<\infty$ and let $\rho_t:=S_t\rho$. Then $\Psi(\rho_t)<\infty$ and $P(\rho_t)<\infty$ for every $t>0$, and we have
\[
\label{eq:cap4:15}
\frac d{dt}\Psi(\rho_t)=-P(\rho_t)\quad\mbox{and}\quad\frac d{dt}P(\rho_t)+2\lambda\,P(\rho_t)\le 0\;.
\]
As a consequence, we have
\[
\label{eq:cap4:16}
\Psi(\rho_t)\le e^{-2\lambda t}\,\Psi(\rho)\;,\quad t\,P(\rho_t)\le (1+2\lambda t)\,e^{-2\lambda t}\,\Psi(\rho)\;,\quad P(\rho_t)\le e^{-2\lambda t}\,P(\rho)
\]
for any $t\ge 0$ and the following \emph{entropy -- entropy production inequality}, or \emph{generalized Poincar\'e inequality}, holds
\[
\label{eq:Cap2:22}
\Psi(\rho)\le \frac 1{2\,\lambda}\,P(\rho)\;,\quad\forall\;\rho\in\mathcal D(\Psi)\quad\mbox{such that}\quad\int_\Rd|x|^2\,\rho\,d\gamma<\infty\;.
\]
\end{theorem}\smallskip
\begin{proof}
It is not restrictive to assume that $\int_\Rd \rho\,d\gamma=1$. We first prove Theorem~\ref {thm:Ent-Prod} for a function $h$ which grows at least linearly at $\infty$, and therefore satisfies $h(r)\ge \sfh\,r$ for some constant $\sfh >0$. The general result follows by writing $h$ as the limit of a decreasing sequence of such concave functions $h_n$, observing that the corresponding actions $\phi_n$ and entropies $\psi_n$ converge increasingly to $\phi$ and $\psi$ respectively.

By \eqref{eq:scheme3:12} we know that $\mathcal H(\rho)$ is finite and therefore we have
\[
\label{eq:Cap2:19}
\int_0^{\infty} \mathcal I(\rho_t)\,dt\le\mathcal H(\rho)<\infty\quad\mbox{and}\quad\int_0^{\infty}\frac{|\dD\,L_\psi(\rho_t)|^2}{\rho_t}\,d\gamma<\infty\;,
\]
where the second estimate follows from \eqref{eq:Cap2:16}.

Applying the chain rule for convex functionals in Wasserstein spaces (see for instance~\cite[p. 233]{ags}, we obtain that the map $t\mapsto \Psi(\rho_t)$ is absolutely continuous and
\[
\label{eq:Cap2:20}
-\frac{d}{dt}\Psi(\rho_t)=\int_\Rd\frac {\dD\,L_\psi(\rho_t)}{\rho_t}\cdot\frac{\dD\rho_t}{\rho_t}\,\rho_t\,d\gamma=P(\rho_t)\;.
\]
By combining Theorems \ref{thm:density_gradient_link} and \ref{thm:action_decay} applied with $\ww_t:=\dD\rho_t$ and differentiating with respect to~$t$, we get that $-\frac d{dt}P(\rho_t)\ge2\lambda\,P(\rho_t)$. All other estimate are easy consequences that have already been established in Section~\ref{Sec:Intro}.
\end{proof}

\section{Contraction of the $h$-Wasserstein distance and KFP as a gradient flow}\label{Sec:Wasserstein}

Consider the space $\mathscr P_{h,\gamma}(\Rd)$ of probability measures at finite $W_{h,\gamma}$ distance from $\gamma$. From \eqref{eq:scheme3:3}, we know that $\gamma$ has finite quadratic moments and, as a consequence of \cite[Theorem 5.9]{MR2435196}, any measure in $\mathscr P_{h,\gamma}(\Rd)$ also has finite quadratic moments. The same result holds for moments of higher order.

\smallskip\begin{theorem}\label{thm:main_contraction}
For every $\sigma,\eta\in\mathscr P_{h,\gamma}(\Rd)$, we have
\[
\label{eq:Cap4:8}
W_{h,\gamma}(\mathcal S_t\sigma\;,\mathcal S_t\eta)\le e^{-\lambda t}\,W_{h,\gamma}(\sigma,\eta)\quad\forall\;t\ge0\;.
\]
\end{theorem}\smallskip
\begin{proof} 
It is a straightforward consequence of Corollary~\ref{cor:action_decay2} and Theorem~\ref{existence:geodesics}.
\end{proof}

\smallskip\begin{theorem}
\label{thm:metric_characterization}
For every $\mu\in\mathscr P_{h,\gamma}(\Rd)$, we have
\begin{equation}
\label{eq:Cap5:1}
\frac 12\frac d{dt}W_{h,\gamma}^2(\mathcal S_t\mu,\sigma)+\frac\lambda2\,W_{h,\gamma}^2(\mathcal S_t\mu,\sigma) + \Psi(\mathcal S_t\mu\,|\,\gamma)\le\Psi(\sigma\,|\,\gamma)\quad\forall\;\sigma\in\mathcal D(\Psi)\;.
\end{equation}
\end{theorem}\smallskip
\begin{proof}
Let us first notice that since $\mathscr P_{h,\gamma}(\Rd)$ is stable under the action of the semigroup $(\mathcal S_t)$, it is sufficient to prove \eqref{eq:Cap5:1} only at $t=0$, under the assumption that~$\mu$ writes as $\mathcal S_\tau\tilde\mu$, for some $\tau>0$. We make the additional assumption on the function~$h$ that there exists some $\mathsf h>0$ for which
\begin{equation}
\label{eq:Cap6:2}
h(r)\geq \mathsf h\,r\quad\forall\;r>0\;.
\end{equation}
This assumption will be removed later in the proof. Let $\varepsilon>0$ fixed and $(\rho^s,\ww^s)\in L^1(\Rd)\times L^1(\Rd,\Rd)$, $s\in[0,1]$, be an admissible curve connecting $\sigma$ to $\mu$ such that
\[
W^2_{h,\gamma}(\mu,\sigma) \leq \mathscr E_\Phi(\rho^s,\ww^s) \leq W^2_{h,\gamma}(\mu,\sigma)+\varepsilon\;,
\]
where $\mathscr E_\Phi(\rho^s,\ww^s):=\int_0^1\Phi(\rho^s,\ww^s)\,ds$. For any $\kappa>0$, we take $\rho^s_\kappa=\rho^s+\kappa\geq\kappa$. Since $h$ is non decreasing, we still have
\begin{equation}
\label{eq:Cap6:3}
\mathscr E_\Phi(\rho^s_\kappa,\ww^s)\leq W^2_{h,\gamma}(\mu,\sigma)+\varepsilon\;.
\end{equation}
Notice that, thanks to \cite[Theorem 5.17]{MR2448650} (also see Theorem~\ref{existence:geodesics}), it is possible to assume that
\begin{equation}
\label{eq:EPhiGeodesics}
\mathscr E_\Phi(\rho^s_\kappa,\ww^s)=\Phi(\rho^s_\kappa,\ww^s)
\end{equation}
is constant with respect to $s\in[0,1]$. For $t>0$, we set
\[
\left\{
\begin{array}{l}
\rho^{s,t}_\kappa = S_{st}\rho^s_\kappa\;,\\[6pt]
\ww^{s,t}_\kappa = \rR_{st}\ww^s-t\,\dD\rho^{s,t}_\kappa\;.
\end{array}
\right.
\]
It is clear that $(\rho^{s,t}_\kappa,\ww^{s,t}_\kappa)$ connects $\sigma+\kappa\gamma$ to $\mathcal S_t(\mu+\kappa\gamma)=\mathcal S_t\mu+\kappa\gamma$. Note that, thanks to the maximum principle, we have $\rho^{s,t}_\kappa\geq \kappa$. We claim that it is admissible. Indeed,
\begin{eqnarray*}
\partial_s\rho^{s,t}_\kappa & = & S_{st}(\partial_s\rho^s_\kappa)+t\,\partial_\tau(S_\tau\rho^s_\kappa)_{|\tau=st}\\ [6pt]
& = & -S_{st}(\tdiv\ww^s)+ \tdiv(t\,\dD\rho^{s,t}_\kappa)\;,
\end{eqnarray*}
since $(\rho^s_\kappa,\ww^s)$ is admissible. Hence,
\[
\partial_s\rho^{s,t}_\kappa = \tdiv\left(-\rR_{st}\ww^s+t\,\dD\rho^{s,t}_\kappa\right)=-\tdiv(\ww^{s,t}_\kappa)\;.
\]
It follows from the definition of $W^2_{h,\gamma}$ that
\[
\label{eq:Cap6:4}
W^2_{h,\gamma}(\mathcal S_t\mu+\kappa\gamma,\sigma+\kappa\gamma) \leq \mathscr E_\Phi(\rho^{s,t}_\kappa,\ww^{s,t}_\kappa)\;,
\]
hence, with \eqref{eq:Cap6:3} and \eqref{eq:EPhiGeodesics}, we obtain
\begin{equation}
\label{eq:Cap6:5}
\frac12\left[W^2_{h,\gamma}(\mathcal S_t\mu+\kappa\gamma,\sigma+\kappa\gamma)-W^2_{h,\gamma}(\mu,\sigma)\right]\leq \frac 12\left[\mathscr E_\Phi(\rho^{s,t}_\kappa,\ww^{s,t}_\kappa)-\mathscr E_\Phi(\rho^s_\kappa,\ww^s)\right]+\frac{\varepsilon}{2}\;.
\end{equation}
By definition of $\mathscr E_\Phi$, we have
\[
\mathscr E_\Phi(\rho^{s,t}_\kappa,\ww^{s,t}_\kappa)=\int_0^1\Phi\left(S_{st}\rho^s_\kappa,\rR_{st}\ww^s-t\,\dD\rho^{s,t}_\kappa\right)ds\;,
\]
where
\begin{multline*}
\Phi\left(S_{st}\rho^s_\kappa,\rR_{st}\ww^s-t\,\dD\rho^{s,t}_\kappa\right)=\int_{\Rd}\frac{\left|\rR_{st}\ww^s - t\,\dD\rho^{s,t}_\kappa\right|^2}{h(\rho^{s,t}_\kappa)}\,d\gamma\\
=\int_{\Rd}\frac{\left|\rR_{st}\ww^s\right|^2}{h(\rho^{s,t}_\kappa)}\,d\gamma-2t\int_{\Rd}\frac{\dD\rho^{s,t}_\kappa\cdot \rR_{st}\ww^s}{h(\rho^{s,t}_\kappa)}\,d\gamma+t^2\int_{\Rd}\frac{\left|\dD\rho^{s,t}_\kappa\right|^2}{h(\rho^{s,t}_\kappa)}\,d\gamma\\
\leq\int_{\Rd}\frac{\left|\rR_{st}\ww^s\right|^2}{h(\rho^{s,t}_\kappa)}\,d\gamma-2t\int_{\Rd}\frac{\dD\rho^{s,t}_\kappa\cdot \ww^{s,t}_\kappa}{h(\rho^{s,t}_\kappa)}\,d\gamma\;,
\end{multline*}
and hence,
\begin{equation}
\label{eq:Cap6:11}
\mathscr E_\Phi(\rho^{s,t}_\kappa,\ww^{s,t}_\kappa) \leq \mathscr E_\Phi(S_{st}\rho^s_\kappa,\rR_{st}\ww^s)-2t\int_0^1\int_{\Rd}\frac{\dD\rho^{s,t}_\kappa\cdot \ww^{s,t}_\kappa}{h(\rho^{s,t}_\kappa)}\,d\gamma\,ds\;.
\end{equation}
\smallskip\begin{lemma}\label{lem:Cap6:chainrule}
If \eqref{eq:Cap6:2} holds, then we have
\begin{equation}
\label{eq:Cap6:9}
\int_0^1\int_{\Rd}\frac{\dD\rho^{s,t}_\kappa\cdot\ww^{s,t}_\kappa}{h(\rho^{s,t}_\kappa)}\,d\gamma\,ds = \Psi(\mathcal S_t\mu+\kappa\gamma\,|\,\gamma) - \Psi(\sigma+\kappa\gamma\,|\,\gamma)\;.
\end{equation}
\end{lemma}\smallskip
\begin{proof}
Recall that $\rho^{s,t}_\kappa=S_{st}\rho^s_\kappa$, with $\rho^s_\kappa\in L^1_\gamma(\Rd)$. Then, acting as in the proof of Proposition~\ref{thm:production}, we get that
\begin{equation}
\label{eq:Cap6:6}
\int_\tau^1\int_{\Rd}\frac{\left|\dD\rho^{s,t}_\kappa\right|^2}{\rho^{s,t}_\kappa}\,d\gamma\,ds\,<\,\infty\quad\forall\;\tau>0\;.
\end{equation}
The assumption \eqref{eq:Cap6:2} on $h$ then leads to
\[
\label{eq:Cap6:7}
\int_\tau^1\int_{\Rd}\frac{\left|
  \dD L_\psi(\rho^{s,t}_\kappa)\right|^2}{\rho^{s,t}_\kappa}\,d\gamma\,ds \leq \frac{1}{\mathsf h^2}\int_\tau^1\int_{\Rd}\frac{\left|\dD\rho^{s,t}_\kappa\right|^2}{\rho^{s,t}_\kappa}\,d\gamma\,ds\,<\,\infty\;.
\]
The next step consists in proving that
\begin{equation}
\label{eq:Cap6:8}
\int_\tau^1\int_{\Rd}\frac{\left|\ww^{s,t}_\kappa\right|^2}{\rho^{s,t}_\kappa}\,d\gamma\,ds\,<\,\infty\;.
\end{equation}
Note that $\rho^{s,t}_\kappa\geq\kappa$ and the concavity of $h$ implies that
\[
h(\rho^{s,t}_\kappa)\leq\frac{h(\kappa)}{\kappa}\,\rho^{s,t}_\kappa\;,
\]
hence
\begin{multline*}
\int_\tau^1\int_{\Rd}\frac{\left|\rR_{st}\ww^s\right|^2}{\rho^{s,t}_\kappa}\,d\gamma\,ds \leq \frac{h(\kappa)}{\kappa}\int_\tau^1\int_{\Rd}\frac{\left|\rR_{st}\ww^s\right|^2}{h(\rho^{s,t}_\kappa)}\,d\gamma\,ds\\
\leq \frac{h(\kappa)}{\kappa}\int_\tau^1\int_{\Rd}\frac{\left|\ww^s\right|^2}{h(\rho^s_\kappa)}\,d\gamma\,ds < \infty
\end{multline*}
since the KFP flow decreases the action. The bound \eqref{eq:Cap6:8} immediately follows from the previous one and \eqref{eq:Cap6:6}. As a consequence, we can apply the chain rule in Wasserstein space, which implies that the function $s\mapsto \Psi(\rho^{s,t}_\kappa)$ is absolutely continuous on $[\tau,1]$ and, for all $s\in[\tau,1]$,
\begin{equation}
\label{eq:Cap6:10}
\frac{d}{ds}\int_{\Rd}\psi(\rho^{s,t}_\kappa)\,d\gamma = \int_{\Rd}
\frac{
 \dD L_\psi(\rho^{s,t}_\kappa)}{\rho^{s,t}_\kappa}\cdot\frac{\ww^{s,t}_\kappa}{\rho^{s,t}_\kappa}\,\rho^{s,t}_\kappa \,d\gamma = \int_{\Rd}\frac{\dD\rho^{s,t}_\kappa\cdot\ww^{s,t}_\kappa}{h(\rho^{s,t}_\kappa)}\,d\gamma\;.
\end{equation}
Integrating \eqref{eq:Cap6:10} on $[\tau,1]$ and letting $\tau$ go to $0$ finally leads to \eqref{eq:Cap6:9}.
\end{proof}

\medskip Let us go back to the proof of Theorem \ref{thm:metric_characterization}. We put \eqref{eq:Cap6:5} and \eqref{eq:Cap6:11} together and obtain
\begin{multline*}
\frac 12\left[W^2_{h,\gamma}(\mathcal S_t\mu+\kappa\gamma,\sigma+\kappa\gamma)-W^2_{h,\gamma}(\mu,\sigma)\right)\\
\leq\frac 12\left[\mathscr E_\Phi(S_{st}\rho^s_\kappa,\rR_{st}\ww^s)-\mathscr E_\Phi(\rho^s_\kappa,\ww^s)\right]
+\,t\left[\Psi(\sigma+\kappa\gamma\,|\,\gamma)-\Psi(\mathcal S_t\mu+\kappa\gamma\,|\,\gamma)\right]+\frac{\varepsilon}{2}\;.
\end{multline*}
We then use the main estimate in Theorem~\ref{thm:action_decay} with $\beta=0$ and \eqref{eq:EPhiGeodesics} to write
\begin{eqnarray*}
\frac 12\left[\mathscr E_\Phi(S_{st}\rho^s_\kappa,\rR_{st}\ww^s)-\mathscr E_\Phi(\rho^s_\kappa,\ww^s)\right] & \leq & -\frac12\,I_\lambda(t)\,\mathscr E_\Phi(\rho^s_\kappa,\ww^s)\\
& \leq & -\frac12\,I_\lambda(t)\,W^2_{h,\gamma}(\mathcal S_t\mu+\kappa\gamma,\sigma+\kappa\gamma)\;,
\end{eqnarray*}
where $I_\lambda(t) := \int_0^1\left(1-e^{-2\lambda ts}\right)ds$. It follows that
\begin{multline*}
\frac 12\left[W^2_{h,\gamma}(\mathcal S_t\mu+\kappa\gamma,\sigma+\kappa\gamma)-W^2_{h,\gamma}(\mu,\sigma)\right] + \frac12\,I_\lambda(t)\,W^2_{h,\gamma}(\mathcal S_t\mu+\kappa\gamma,\sigma+\kappa\gamma)\\
\label{eq:Cap6:13} \leq
t\left[\Psi(\sigma+\kappa\gamma\,|\,\gamma)-\Psi(\mathcal S_t \mu+\kappa\gamma\,|\,\gamma)\right]+\frac{\varepsilon}{2}\;.
\end{multline*}
If we first let $\varepsilon$ and then $\kappa$ go to $0$ in the above estimate, we get that
\begin{equation}
\label{eq:Cap6:14}
\frac 12\left[W^2_{h,\gamma}(\mathcal S_t\mu,\sigma)-W^2_{h,\gamma}(\mu,\sigma)\right]+\frac12\,I_\lambda(t)\,W^2_{h,\gamma}(\mathcal S_t\mu,\sigma) \leq t\left[\Psi(\sigma\,|\,\gamma)-\Psi(\mathcal S_t\mu\,|\,\gamma)\right]
\end{equation}
as soon as $h$ satisfies the assumption \eqref{eq:Cap6:2}. Now, any concave and non decreasing function $h$ can be decreasingly approched by a sequence $(h_n)$ satisfying \eqref{eq:Cap6:2}, and the corresponding entropies converge increasingly. Then, with Theorem~\ref{thm:2}, Inequality~\eqref{eq:Cap6:14} turns out to be valid for any general $h$. To complete the proof of Theorem~\ref{thm:metric_characterization}, it just remains to divide \eqref{eq:Cap6:14} by $t$ and let $t$ go to $0$.
\end{proof}

\medskip\noindent {\bf Acknowledgements.} This work has been partially
supported by the projects CBDif, EVOL and OTARIE of the French National
Research Agency (ANR) and by a PRIN08-grant from MIUR for the project
\emph{Optimal transport theory, geometric and functional inequalities,
 and applications.}
\par\smallskip\noindent{\small\copyright\,2011 by the authors. This paper may be reproduced, in its entirety, for non-commercial purposes.}

\def\cprime{$'$}


\begin{thebibliography}{10}

\bibitem{ambrosio-Fusco-Pallara00}
{\sc L.~Ambrosio, N.~Fusco, and D.~Pallara}, {\em Functions of bounded
  variation and free discontinuity problems}, Oxford Mathematical Monographs,
  Clarendon Press, Oxford, 2000.

\bibitem{Ambrosio-Gigli-Savare05}
{\sc L.~Ambrosio, N.~Gigli, and G.~Savar{\'e}}, {\em Gradient flows in metric
  spaces and in the space of probability measures}, Lectures in Mathematics ETH
  Z\"urich, Birkh\"auser Verlag, Basel, 2005.

\bibitem{ags}
\leavevmode\vrule height 2pt depth -1.6pt width 23pt, {\em Gradient flows in
  metric spaces and in the space of probability measures}, Lectures in
  Mathematics ETH Z\"urich, Birkh\"auser Verlag, Basel, second~ed., 2008.

\bibitem{Ambrosio-Savare06}
{\sc L.~Ambrosio and G.~Savar\'e}, {\em Gradient flows of probability
  measures}, in Handbook of Evolution Equations (III), Elsevier, 2006.

\bibitem{Ambrosio-Savare-Zambotti07}
{\sc L.~Ambrosio, G.~Savar{\'e}, and L.~Zambotti}, {\em Existence and stability
  for {F}okker-{P}lanck equations with log-concave reference measure}, Probab.
  Theory Related Fields, 145 (2009), pp.~517--564.

\bibitem{MR2375056}
{\sc A.~Arnold, J.-P. Bartier, and J.~Dolbeault}, {\em Interpolation between
  logarithmic {S}obolev and {P}oincar\'e inequalities}, Commun. Math. Sci., 5
  (2007), pp.~971--979.

\bibitem{MR2152502}
{\sc A.~Arnold and J.~Dolbeault}, {\em Refined convex {S}obolev inequalities},
  J. Funct. Anal., 225 (2005), pp.~337--351.

\bibitem{Arnold-Markowich-Toscani-Unterreiter01}
{\sc A.~Arnold, P.~Markowich, G.~Toscani, and A.~Unterreiter}, {\em On convex
  {S}obolev inequalities and the rate of convergence to equilibrium for
  {F}okker-{P}lanck type equations}, Comm. Partial Differential Equations, 26
  (2001), pp.~43--100.

\bibitem{MR889476}
{\sc D.~Bakry and M.~{\'E}mery}, {\em Diffusions hypercontractives}, in
  S\'eminaire de probabilit\'es, {XIX}, 1983/84, vol.~1123 of Lecture Notes in
  Math., Springer, Berlin, 1985, pp.~177--206.

\bibitem{MR2317340}
{\sc D.~Bakry, M.~Ledoux, and F.-Y. Wang}, {\em Perturbations of functional
  inequalities using growth conditions}, J. Math. Pures Appl. (9), 87 (2007),
  pp.~394--407.

\bibitem{MR2320410}
{\sc F.~Barthe, P.~Cattiaux, and C.~Roberto}, {\em Interpolated inequalities
  between exponential and {G}aussian, {O}rlicz hypercontractivity and
  isoperimetry}, Rev. Mat. Iberoam., 22 (2006), pp.~993--1067.

\bibitem{MR2346509}
{\sc F.~Barthe, P.~Cattiaux, and C.~Roberto}, {\em Isoperimetry between
  exponential and {G}aussian}, Electron. J. Probab., 12 (2007), pp.~no. 44,
  1212--1237.

\bibitem{Beckner89}
{\sc W.~Beckner}, {\em A generalized {P}oincar\'e inequality for {G}aussian
  measures}, Proc. Amer. Math. Soc., 105 (1989), pp.~397--400.

\bibitem{Benamou-Brenier00}
{\sc J.-D. Benamou and Y.~Brenier}, {\em A computational fluid mechanics
  solution to the {M}onge-{K}antorovich mass transfer problem}, Numer. Math.,
  84 (2000), pp.~375--393.

\bibitem{MR2609029}
{\sc F.~Bolley and I.~Gentil}, {\em Phi-entropy inequalities for diffusion
  semigroups}, J. Math. Pures Appl. (9), 93 (2010), pp.~449--473.

\bibitem{MR2081075}
{\sc D.~Chafa{\"{\i}}}, {\em Entropies, convexity, and functional inequalities:
  on {$\Phi$}-entropies and {$\Phi$}-{S}obolev inequalities}, J. Math. Kyoto
  Univ., 44 (2004), pp.~325--363.

\bibitem{MR2452882}
{\sc S.~Daneri and G.~Savar{\'e}}, {\em Eulerian calculus for the displacement
  convexity in the {W}asserstein distance}, SIAM J. Math. Anal., 40 (2008),
  pp.~1104--1122.

\bibitem{MR2273884}
{\sc P.~Deng and F.~Wang}, {\em Beckner inequality on finite- and
  infinite-dimensional manifolds}, Chinese Ann. Math. Ser. B, 27 (2006),
  pp.~581--594.

\bibitem{MR2366398}
{\sc J.~Dolbeault, I.~Gentil, A.~Guillin, and F.-Y. Wang}, {\em {$L\sp
  q$}-functional inequalities and weighted porous media equations}, Potential
  Anal., 28 (2008), pp.~35--59.

\bibitem{MR2435196}
{\sc J.~Dolbeault, B.~Nazaret, and G.~Savar{\'e}}, {\em On the {B}akry-{E}mery
  criterion for linear diffusions and weighted porous media equations}, Commun.
  Math. Sci., 6 (2008), pp.~477--494.

\bibitem{MR2448650}
{\sc J.~Dolbeault, B.~Nazaret, and G.~Savar{\'e}}, {\em A new class of
  transport distances between measures}, Calc. Var. Partial Differential
  Equations, 34 (2009), pp.~193--231.

\bibitem{Gianazza-Savare-Toscani09}
{\sc U.~Gianazza, G.~Savar{\'e}, and G.~Toscani}, {\em The {W}asserstein
  gradient flow of the {F}isher information and the quantum drift-diffusion
  equation}, Arch. Ration. Mech. Anal., 194 (2009), pp.~133--220.

\bibitem{john}
{\sc F.~John}, {\em Hyperbolic and parabolic equations}, in Partial
  {D}ifferential {E}quations ({P}roc. {S}ummer {S}eminar, {B}oulder, {C}ol.,
  1957), Interscience, New York, 1964, pp.~1--129.

\bibitem{Jordan-Kinderlehrer-Otto98}
{\sc R.~Jordan, D.~Kinderlehrer, and F.~Otto}, {\em The variational formulation
  of the {F}okker-{P}lanck equation}, SIAM J. Math. Anal., 29 (1998),
  pp.~1--17.

\bibitem{MR1796718}
{\sc R.~Lata{\l}a and K.~Oleszkiewicz}, {\em Between {S}obolev and
  {P}oincar\'e}, in Geometric aspects of functional analysis, vol.~1745 of
  Lecture Notes in Math., Springer, Berlin, 2000, pp.~147--168.

\bibitem{MR2192294}
{\sc F.~Otto and M.~Westdickenberg}, {\em Eulerian calculus for the contraction
  in the {W}asserstein distance}, SIAM J. Math. Anal., 37 (2005),
  pp.~1227--1255.

\bibitem{MR2459454}
{\sc C.~Villani}, {\em Optimal transport. Old and new}, vol.~338 of Grundlehren
  der Mathematischen Wissenschaften [Fundamental Principles of Mathematical
  Sciences], Springer-Verlag, Berlin, 2009.

\end{thebibliography}
\end{document}